\documentclass[11pt]{article}
\usepackage{amssymb}
 \usepackage{times}
\usepackage{amscd}

\begin{document}
\newtheorem{theorem}{Theorem}[section]
\newtheorem{lemma}[theorem]{Lemma}
\newtheorem{corollary}[theorem]{Corollary}
\newtheorem{definition}[theorem]{Definition}
\newtheorem{proposition}[theorem]{Proposition}
\newtheorem{defprop}[theorem]{Definition-Proposition}
\newtheorem{example}[theorem]{Example}
\newtheorem{remark}[theorem]{Remark}
\newcommand{\Proof}{\noindent{\bf Proof:} } 
\catcode`\@=11
\@addtoreset{equation}{section}
\catcode`\@=12
\renewcommand{\theequation}{\arabic{section}.\arabic{equation}}
\def\lu{\rightharpoonup}
\def\Infl{{\rm Infl}}
\def\sqr#1#2{{\vcenter{\vbox{\hrule height.#2pt\hbox{\vrule width.#2pt
  height#1pt \kern#1pt \vrule width.#2pt}\hrule height.#2pt}}}}
\def\square{\mathchoice\sqr64\sqr64\sqr{2.1}3\sqr{1.5}3}
\def\argh{\rightharpoonup} 
\def\End{{\rm End}}
\def\m#1{m_{(#1)}}
\def\g#1{g_{(#1)}}
\def\f{{\bar{f}}}
\def\h#1{h_{(#1)}}
\def\ad{{\rm ad}}
\def\k#1{k_{(#1)}}
\def\R#1{R_{(#1)}}
\def\M{{\cal M}} 
\def\Z{{\mathbb Z}}
\def\Z#1{{\mathbb Z}_#1}
\def\F{{\cal F}}
\def\mm#1{m_{(#1*)}}
\def\a#1{a_{(#1)}}
\def\aa#1{a_{(#1*)}}
\def\b#1{b_{(#1)}}
\def\bb#1{b_{(#1*)}}
\def\l#1{l_{(#1)}}
\def\n#1{n_{(#1)}}
\def\nn#1{n_{(#1*)}}
\def\totimes{\tilde{\otimes}}
\def\c#1{c_{(#1)}}
\def\cc#1{c_{(#1*)}}
\def\id{{\rm id}}
\def\ydh{{\cal YD}(H)}
\def\yddh{{\cal YD}(H^*)}
\def\lefth{_{H}{\cal M}}
\def\lefthd{_{H^*}{\cal M}}
\def\righth{{\cal M}^H}
\def\righthd{{\cal M}^{H^*}}
\def\C{{\cal C}}
\def\D{{\cal D}}
\def\rsb{r_{s,\beta}}
\def\rsg{r_{s,\gamma}}
\def\km{k^{*}}
\def\pf{{\bf Proof:} }
\def\U{\langle u\rangle}
\def\H{{\cal H}}
\def\CoInn{{\rm CoInn}}
\def\CoInt{{\rm CoInt}}
\def\Zen{Z_L(E(n))}
\def\sym{Sym_{M,n,r}(k)}
\def\Hom{{\rm Hom}}
\def\Ker{{\rm Ker}}
\title{The Brauer group of modified supergroup
  algebras} 
\author{Giovanna Carnovale\\
Dipartimento di Matematica Pura ed Applicata\\
Universit\`a di Padova\\
via Belzoni 7\\
I-35131 Padova, Italy\\
carnoval@math.unipd.it} 
\date{}
\maketitle
\begin{abstract}
The computation of the Brauer group $BM$ of modified supergroup
algebras is perfomed, yielding, in particular, the computation of the
Brauer group of all finite-dimensional triangular Hopf algebras when
the base field is algebraically closed and of characteristic zero. The
results are compared with the computation of lazy cohomology and with Yinhuo Zhang's exact sequence. 
As an example, we compute explicitely the Brauer group and
lazy cohomology for modified supergroup algebras with (extensions of) 
Weyl groups of irreducible root systems as a group datum and their standard
representation as a representation datum.
\end{abstract} 

 
\section*{Introduction}

Finite-dimensional triangular Hopf algebras over an algebraically
closed base field of characteristic zero have been completely
described in \cite{chev1},\cite{EG3},\cite{EG2},\cite{EG1}. They
can all be reduced, by Drinfeld twists (the construction dual to
a cocycle twist), to a
particular class of pointed Hopf algebras,  the so-called modified
supergroup algebras. Such algebras are directly constructed starting
from a finite group $G$, a central involution $u$ of $G$, and a
representation $V$ of $G$ on which $u$ acts as $-1$.  
The $R$-matrix can always be reduced to $R_u=\frac{1}{2}(1\otimes
1+1\otimes u+u\otimes 1-u\otimes u)$. 

In \cite{CVZ} and \cite{CVZ2} several group invariants 
for finite-dimensional Hopf
algebras have been defined. In particular, the Brauer group
$BM(k,H,R)$ of a quasitriangular Hopf algebra $(H,R)$ has been
defined as the Brauer group of the braided monoidal category of left
$H$-modules. 
It has been shown in \cite{gio} as a
consequence of results in \cite{VZ} that this Brauer group is invariant under
Drinfeld twists of $(H, R)$  and this fact can be used in order to replace $(H,R)$ with a new pair $(H', R')$
 that is easier to deal with. In particular, this property has been
 used to replace the $R$-matrix $R$ with a simpler one.  
In the last few years several explicit
examples have been computed for triangular or quasitriangular 
Hopf algebras. All of these Hopf algebras 
admitted several distinct quasitriangular
structures and most of them admitted at least one
triangular structure: they were in fact modified supergroup algebras. 
The first explicit computation was performed in \cite{VOZ3} where the
Brauer group of 
Sweedler Hopf algebra $H_4$ is determined when the triangular
structure is $R_u$. This computation was generalized to all possible
(quasi)triangular
structures of $H_4$ in \cite{gio}.  In terms of modified supergroup
algebras, $H_4$ corresponds to the data $G={\mathbb Z}_2$ and $V$ 
its non-trivial irreducible representation.
In \cite{hnu} the computation of $BM$ is generalized to the case in which $G$
is the cyclic group $\langle g\rangle$ of order $2\nu$ for $\nu$ an
odd integer and $V$ is 
the irreducible representation on which $g$ again acts as $-1$. In
 \cite{diedrale} the first case of a non-abelian group was considered
 and in \cite{GioJuan3} the case in which $G={\mathbb Z}_2$ and $V$ is
 given by $n$-copies of the non-trivial irreducible representation of
 $G$ was treated, for all triangular structures. In all those cases 
important roles were played by the
 classical Brauer group, by the Brauer-Wall group, and by lazy cohomology as
 defined  
in \cite{scha} (with the name central cohomology) and studied in
\cite{lazy} and \cite{Juan-Florin}.
It is therefore natural to generalize the approach used in the
aforementioned papers in order to deal
with all $G$ and $V$ at the same time, for $R$ triangular, obtaining
in particular
the computation of the Brauer group of all finite-dimensional
triangular Hopf algebras over an algebraically closed field of
characteristic zero. This construction should also give
an indication on how to compute the Brauer group of quasitriangular
Hopf algebras admitting also a triangular structure. 

We list here the content of the paper. 
First, the case in which $V=0$ is dealt with. In this case the
Brauer group of the group algebra $k[G]$ with $R$-matrix $R_u$ 
is
\begin{itemize}
\item  isomorphic to $Br(k)\times H^2(G,k^\cdot)$ if $u=1$;
\item a central extension of $Br(k)$ by  $H^2_\sharp(G,
k^\cdot)$ if $u\neq1$ and $U=\langle u\rangle$ is not a direct summand
of $G$. Here, for a group $K$, the symbol 
$H^2_\sharp(K, k^\cdot)$ denotes an abelian group coinciding with $H^2(K,
k^\cdot)$ as a set but with a modified group structure given
in Section \ref{sharp}. Conditions on $K$ or $k$ for  $H^2_\sharp(K,
k^\cdot)$ to be isomorphic to  $H^2(K,
k^\cdot)$ are given in Proposition \ref{coincide}. For instance, the two
groups coincide when $-1$ is a square in $k$. 
\item a central extension of $Br(k)$ by
$Q(k,G)$, which is itself a central extension of 
$H^2_\sharp(G,k^\cdot)$ by ${\mathbb Z}_2$, if  $u\neq1$ and $U=\langle
  u\rangle$ is a direct summand 
of $G$.
\end{itemize}  

Thanks to the results in \cite{EG1} we obtain:

\noindent{\bf Theorem }{\em Let $(H,R)$ be a semi-simple and
cosemi-simple triangular 
Hopf algebra over $k=\bar{k}$ with ${\rm gcd}({\rm char}(k),\dim H)=1$. Then 
$$
BM(k, H, R)=\cases{H^2(G, k^\cdot)\times{\mathbb Z}_2& if
$u\neq1$ and $G\cong \U\times G/\U$;\cr
H^2(G, k^\cdot)& otherwise\cr} 
$$
where $G,\,u$ and $V=0$ are the corresponding data for $(H, R)$.}

\medskip 

The more general case in which $V\neq0$ is dealt with by 
showing that the Brauer group in this case is the direct
product of the Brauer group of $k[G]$ and the group of $G$-invariant
symmetric bilinear forms on $V^*$. In this case $R=R_u$ with $u\neq1$. By \cite[Proposition 2.6]{EG2} we obtain:

\medskip

\noindent{\bf Theorem }{\em  Let $(H, R)$ be a finite-dimensional
triangular Hopf 
algebra over $k={\bar k}$ with ${\rm char}(k)=0$. Then  
$$BM(k, H, R)=\left\{
\begin{array}{lc}{\mathbb Z}_2\times H^2(G,
    k^\cdot)\times S^2(V^*)^G& 
\mbox{if $u\neq1$ and } G\cong
    U\times G/U\\
&\\
H^2(G, k^\cdot)\times S^2(V^*)^G &\mbox{ otherwise}\\
\end{array}\right.
$$ 
where $G,\,u,\,V, U=\langle u\rangle$ are the corresponding data for $(H, R)$.}

\medskip

As in the classical case, a link is expected between a second
``cohomology group'' for the Hopf algebra $H$ and its full Brauer
group, that is, the Brauer group of the Drinfeld double $D(H)$ of $H$.  
It is also expected, and evidence was given by all computations made so
far, that this cohomology group could be given by the second lazy
cohomology group $H^2_L(H)$. We adapt the analysis we made for the Brauer group
in order to compute $H^2_L(H)$ for modified supergroup algebras and
compare it to $BM(k, H, R_u)$. If the representation datum $V$ is trivial
then $H=k[G]$ and  
$H^2_L(H)$ is the usual cohomology in degree two of the group datum $G$. 
If the representation datum $V$ is non-trivial
then  
$H^2_L(H)$ is isomorphic to the direct product of $S^2(V^*)^G$ and 
the second cohomology group of $G/U$. Thus, the linear part of lazy
cohomology is a direct summand of $BM(k, H, R_u)$, while for the group
cohomology component the situation is more involved. This phenomenon
is not completely new because already when $H=k[{\mathbb Z}_2]$ the
cohomology of ${\mathbb Z}_2$ with trivial coefficients can be seen as
a subquotient but not as a subgroup of $BM(k, H, R_u)\cong BW(k)$, the
Brauer-Wall group of the base field $k$.   
These relations however cannot be generalized directly to all
finite-dimensional triangular 
Hopf algebras because $H^2_L(H)$ does not seem to be invariant under
Drinfeld twists but only under cocycle twists. 
However, an analysis of the cases in the literature where $H_L^2(H)$ 
has been computed shows that, for instance, for all finite-dimensional 
triangular Hopf algebras $H$ over $k=\bar k$  of 
characteristic zero, for which the $G$-action is faithful,
 $H^2_L(H^*)$ is a direct summand of $BM(k, H, R)$.

We relate the computation of $BM$ with the exact sequence
in \cite{sequence} which generalizes sequences in \cite{beattie},
\cite{CGO}, \cite{ulbrich}. It involves the Brauer group
$BC(k,H,R)\cong BM(k, H^*, R)$ of a dual quasitriangular
 Hopf algebra $(H, R)$ and a group of quantum
 commutative biGalois objects. The results contained in this paper
 should facilitate the interpretation of the results in
 \cite{sequence} and their comparison should serve as an indication
 for handling the quasitriangular case.

In order to illustrate the preceding results, we deal
explicitely with the case of the modified supergroup algebras over
${\mathbb C}$ with (an extension of) the Weyl group corresponding
to an irreducible root system $\Phi$ as a group datum and its standard
representation as a representation datum. Here, explicit computations
are simpler because the representation $V$ of $G$ is faithful and
irreducible, the classical Schur multiplier of all such groups is
known, and, due to the particular presentation of $G$ as a Coxeter group,
most of the different cases occurring in the general
case coincide. 

\section{Preliminaries}

In this Section we shall introduce the basic notions that we shall use
in this paper. 
Unless otherwise stated, $H$ shall denote a
finite-dimensional Hopf algebra over the base field $k$ with product
$m$,  coproduct
$\Delta(h)=\sum\h1\otimes\h2$ and antipode $S$. We shall always assume
that ${\rm gcd}(\dim H, {\rm char}(k))=1$. Unadorned tensor product will
be intended to be over $k$. The symbol $u$ shall always denote a 
central element of a group $G$ for which $u^2=1$ and $U$ shall denote
the subgroup of $G$ generated by $u$. 

\subsection{Modified supergroup algebras}\label{supergroup}

In this subsection we shall recall the notion of a modified supergroup
algebra and its importance within the theory of finite-dimensional
triangular Hopf algebras.

\begin{definition}
A finite dimensional triangular Hopf algebra $H$ with $R$-matrix $R$ 
 is called {\em a
modified supergroup algebra} if there exist:
\begin{itemize}
\item  a finite group $G$, 
\item a central element $u$ of $G$ with $u^2=1$, 
\item a linear representation of $G$ on a finite-dimensional vector space $V$
on which $u$ acts as $-1$,
\end{itemize}  
such that: $H\cong k[G]\ltimes\wedge V$ as an algebra; the
elements in $G$ are grouplike; the elements in $V$ are
$(u,1)$-primitive; $R=R_u=\frac{1}{2}(1\otimes1+1\otimes u+u\otimes
1-u\otimes u)$.      
\end{definition} 

This name is motivated by the fact that such an Hopf algebra is
obtained from a cocommutative Hopf superalgebra through a suitable
modification (see \cite{mabook} where it is called bozonisation, or
\cite{EG3}, \cite{EG2}) and  by the description, in
\cite[Theorem 3.3]{kos}, of all
finite-dimensional cocommutative Hopf superalgebras over ${\mathbb C}$.

Modified supergroup algebras are the model 
for finite-dimensional triangular Hopf algebras, at least over ``good
fields''.  
%
\begin{theorem}\label{classi}(\cite[Theorem 5.1.1]{chev1}, \cite[Theorem
  4.3]{EG1}) Every finite-dimensional
triangular Hopf algebra over an algebraically closed field of
  characteristic zero is the Drinfeld twist of a modified
  supergroup algebra.
\end{theorem}

For the well-known notion of a Drinfeld twist the reader is referred to
\cite[Section 2.4]{chev1}. If $H$ satisfies the conditions of Theorem
\ref{classi} then it is the twist of $k[G]\ltimes\wedge V$ for some
$G$ with central $u$, such that $u^2=1$ and for some representation $V$ on
which $V$ acts as $-1$. We
shall call $G$ the group datum and $V$ the representation datum
corresponding to $H$. 

\begin{example}\label{enne}{\rm If $G={\mathbb Z}_2$ then $k[G]\ltimes\wedge V$
  is isomorphic to the Hopf algebra usually denoted by $E(n)$ with
  $n=\dim(V)$. We recall that $E(n)$ is generated by
$c$ and $x_i$ for
$i=1,\,\ldots,\,n$ with relations
$c^2=1,\ cx_i+x_ic=0,\  x_ix_j+x_jx_i=0,\;\forall i=1,\,\cdots,\,n,$
coproduct 
$\Delta(c)=c\otimes c,\  \Delta(x_i)=1\otimes x_i+x_i\otimes c$
and antipode
$S(c)=c,\quad S(x_j)=cx_j$. We shall use the isomorphism given by $c\mapsto u$
and $x_i\mapsto uv_i$, where $v_1,\ldots,\,v_n$ is a basis of $V$.  
It is well-known that $E(n)$ has a family of triangular $R$-matrices  
parametrized by $S^2(V^*)$, that is, by symmetric $n\times n$ matrices
with coefficients in $k$
(see \cite[Proposition 2.1]{GioJuan3} and \cite{PVO1}). The family is given
as follows. For the symmetric matrix
$A=(a_{ij})$ and for the $s$-tuples $P,\,F$ of increasing elements in
$\{1,\ldots,\,n\}$ we define $|P|=|F|=s$ and $v_P$ as the product of the
$v_j$'s whose index belongs to $P$, taken in increasing order. A map
$\eta$ from the elements of the $s$-uple
$P$ to the elements of the $s$-uple $F$ determines an element of the
symmetric group
$S_s$ which we identify with $\eta$. We denote then by $sign(\eta)$
the sign of $\eta$. If $P$ is empty, i.e., if
$s=0$ 
we take $F$ to be empty and $\eta$ to have sign equal to $1$. Finally
$a_{P,\eta(F)}$ denotes the product $a_{p_1,f_{\eta(1)}}\cdots\,
a_{p_s,f_{\eta(s)}}$. If $P$ is  empty we define
$a_{P,\eta(F)}:=1$. Then the 
$R$-matrix corresponding to $A$ is
$$R_A={1\over2}\sum_P(-1)^{{|P|(|P|-1)}\over2}\sum_{F,\; |F|=|P|,\; \eta\in
S_{|P|}}sign(\eta)a_{P,\eta(F)}\bigl(v_P\otimes
v_F$$
$$+uv_{P}\otimes v_F+(-1)^{|P|}v_P\otimes
uv_F-(-1)^{|P|}v_{P}\otimes uv_F\bigr).$$

It is also well-known that there exists a Hopf algebra isomorphism 
$\phi\colon E(n)\to
E(n)^*$  given by $\phi(1)=\varepsilon$, $\phi(v_j)=v_j^*-(uv_j)^*$ and
$\phi(u)=1^*-u^*$. This is not the case for a general modified
supergroup algebra. By self-duality, $E(n)$ has a family of
dual triangular structures parametrized by symmetric matrices. The family is
given by:
$$r_A=\sum_P (-1)^{{|P|(|P|-1)}\over2}\sum_{F, |F|=|P|, \eta\in
S_{|P|}}sign(\eta)a_{P,\eta(F)}\bigl((v_P)^*\otimes
(v_F)^*$$
$$+(v_{P})^*\otimes (uv_F)^*+(-1)^{|P|}(uv_P)^*\otimes
(v_F)^*-(-1)^{|P|}(uv_{P})^*\otimes (uv_F)^*\bigr).$$
The pair $(E(n), R_A)$ is twist equivalent to the pair $(E(n),
R_0)=(E(n), R_u)$ for every $A$.}
\end{example}

\subsection{The Brauer group $BM(k,H,R)$}\label{definizioni} 

In this Section we shall recall the construction of the Brauer group
of the category $_H{\cal M}$ of left modules for a quasitriangular Hopf
algebra (see \cite{CVZ}, \cite{par}, \cite{VZ} for further details). 
Let $(H, R)$ be a finite-dimensional quasitriangular Hopf algebra over the
base field $k$. Let us write $R=\sum \R1\otimes\R2$. We recall that
$_H{\cal M}$ is a braided 
 monoidal category with braiding given, for all $H$-modules $M,N$ and for all
 $m\in M$ and $n\in N$, by  $$\psi_{MN}:M \otimes N
\rightarrow N \otimes M, \ m \otimes n \mapsto \sum (\R2\cdot
n) \otimes (\R1 \cdot m)$$ for all $m \in M,n \in N.$ 

Given
two $H$-module algebras $A$ and $B$ their {\it braided product} 
$A \sharp B$  is the $H$-module algebra with $A\otimes B$ as underlying
$H$-module and
multiplication given by $$(a \sharp x)(b \sharp y)=a\psi_{BA}(x
\otimes b)y=
\sum a(\R2\cdot b)\sharp (\R1 \cdot x)y,$$ for all $a,b
\in A,x,y \in B.$ Given a $H$-module algebra $A$ its $H$-opposite
algebra $\overline{A}$ is the $H$-module algebra  equal to $A$ as an
$H$-module and with 
multiplication given by $ab=m_A\psi_{AA}(a\otimes b)=\sum (\R2\cdot b)(\R1 \cdot a)$ 
for all $a,b \in A$. 

The endomorphism algebra
$\End(M)$ of a finite-dimensional $H$-module $M$  is an $H$-module
algebra with 
action
\begin{equation}\label{endo}(h \cdot f)(m)=\sum \h1\cdot f(S(\h2)
  \cdot m).\end{equation} 
Its opposite algebra $\End(M)^{op}$ is also a left $H$-module algebra
once it is equipped with the $H$-action:
$$(h \cdot f)(m)=\sum h_{(2)} \cdot f(S^{-1}(h_{(1)}) \cdot m).$$

The following two maps are $H$-module algebra maps:
$$\begin{array}{ll}
F_1: A \sharp \overline{A} \rightarrow \End(A),\ F_1(a \sharp \bar{b})(c)=\sum
a(\R2 \cdot c)(\R1 \cdot b), \smallskip \\
F_2: \overline{A} \sharp A \rightarrow \End(A)^{op},\ F_2(\bar{a} \sharp b
)(c)=\sum (R_{(2)} \cdot a)(R_{(1)} \cdot c)b.
\end{array}$$

We shall say that a finite-dimensional $H$-module algebra $A$ is 
$H$-Azumaya (or $(H, R)$-Azumaya if more $R$-matrices are involved) 
if $F_1$ and $F_2$ are isomorphisms.

The set of isomorphism classes of $H$-Azumaya
algebras is equipped with the equivalence relation: ``$A \sim B$ 
if there exist finite
dimensional $H$-modules $M,N$ such that $A \sharp \End(M) \cong B \sharp
\End(N)$ as $H$-module algebras''. The quotient set $BM(k,H,R)$ inherits
a natural group structure where the multiplication is induced by
$\sharp$, the inverse of a class represented by $A$ is represented by
$\overline{A}$ and the neutral element is represented by $\End(M)$ for
a finite dimensional $H$-module $M$. If $(H, R)$ is triangular then
$BM(k,H, R)$ is abelian. 

\smallskip

The Brauer group $BM(k, H, R)$ embeds naturally into 
$BM(k, D(H), {\cal R})$, that is, into the Brauer group of the
Drinfeld double of $H$. This group is also 
called the full Brauer group of $H$ and it is usually denoted by
$BQ(k, H)$.  

For a dual quasitriangular Hopf algebra $(H',r)$ we can perform the
construction dual to the construction of  
$BM$ obtaining the Brauer group $BC(k, H',r)$ of 
the category ${\cal
M}^H$ of right $H$-comodules. In particular, if $(H, R)$ is a
quasitriangular Hopf algebra then $BM(k, H, R)\cong BC(k,
H^*,R)$. 

\smallskip

We shall often make use of the following result, which is a consequence
of the functoriality of the Brauer group (see \cite{par}, \cite{VZ}):
\begin{theorem}(\cite[Dual of Proposition 3.1]{gio}) The Brauer group
  $BM(k, H, R)$ is 
invariant under Drinfeld twists of $(H, R)$.
\end{theorem}

We end this Section recalling a few notions on the second lazy
cohomology group of a Hopf algebra, as defined in \cite{scha} and
studied in \cite{lazy} and \cite{Juan-Florin}.

We recall that a left 2-cocycle for a Hopf algebra $H$  is a
normalized convolution invertible element $\sigma$ in $(H\otimes H)^*$ 
such that
$$\sum\sigma(\a1,\b1)\sigma(\a2\b2,c)=\sum\sigma(\b1,\c1)\sigma(a,\b2\c2)$$
for every $a,b,c\in H$ and that a right 2-cocycle is a
normalized convolution invertible element in $(H\otimes H)^*$ 
such that
$$\sum\sigma(\a1\b1,c)\sigma(\a2,\b2)=\sum\sigma(a,\b1\c1)\sigma(\b2,\c2)$$
for every $a,b,c\in H$. 
A lazy cocycle is a left 2-cocycle such that
$$\sum\sigma(\a1,\b1)\a2\b2=\sum\sigma(\a2,\b2)\a1\b1$$ for every $a,b\in
H$. 

It is not hard to verify that lazy cocycles are also right cocycles
and that the set of 
lazy cocycles forms a group under convolution which we shall denote by
$Z^2_L(H)$ (\cite[Page 227]{chen}, \cite[Lemma 1.2]{lazy}).  

A left 2-cocycle $\sigma$ in $(H\otimes H)^*$ is
called a left coboundary if
$$\sigma(a,b)=\partial(\gamma)(a,b)=\sum
\gamma(\a1)\gamma(\b1)\gamma^{-1}(\a2\b2)$$ 
for some convolution invertible $\gamma\in H^*$. The set of left
coboundaries will be denoted by $B^2(H)$. Similarly, a right cocycle $\omega$
is called a right coboundary if $$\omega(a,b)=\sum
\gamma(\a1\b1)\gamma^{-1}(\a2)\gamma^{-1}(\b2)$$ for some invertible
$\gamma\in H^*$. 
If $\gamma$ is lazy,
that is, if $\gamma$ is central in $H^*$, then we shall say that
$\sigma=\partial(\gamma)$ is a lazy coboundary. The set of lazy coboundaries forms a
central subgroup of $Z^2_L(H)$ which we shall denote by
$B^2_L(H)$. The factor group $Z^2_L(H)/B^2_L(H)$ is called the second
lazy cohomology group of $H$ and it is denoted by $H^2_L(H)$. 

One should observe that $B^2(H)\cap Z^2_L(H)\neq B^2_L(H)$ in
general. The factor group $B^2(H)\cap Z^2_L(H)/B^2_L(H)$ shall play a
role in the following sections.   

A Hopf algebra automorphism $\psi$ is called {\em coinner} if it is of the
form
$\psi(a)=\ad(\gamma):=\sum \gamma^{-1}(\a1)\a2\gamma(\a3)$ for some
grouplike element $\gamma\in H^*$. The group of coinner automorphism
of $H$ is denoted by $\CoInn(H)$. 

A Hopf algebra automorphism $\psi$ is called {\em cointernal} if it is of the
form
$\psi(a)=\ad(\gamma):=\sum\gamma^{-1}(\a1)\a2\gamma(\a3)$ for some
convolution invertible element $\gamma\in H^*$. 
The group of cointernal automorphism
of $H$ is denoted by $\CoInt(H)$. 
It is proved in \cite[Lemma 1.11]{lazy} that
$\ad(\gamma)$ is a Hopf algebra automorphism if and only if
$\partial(\gamma)\in Z^2_L(H)$.

\begin{lemma}\label{cappa}The groups $B^2(H)\cap Z^2_L(H)/B^2_L(H)$ and
$\CoInt(H)/\CoInn(H)$ are isomorphic. 
\end{lemma}
\pf The proof follows from the proof of \cite[Lemma
1.12]{lazy}. Indeed, the assignment $\partial(\gamma)\mapsto \ad(\gamma)\circ
\CoInn(H)$ is a well-defined surjective group morphism $B^2(H)\cap
Z^2_L(H)\to\CoInt(H)/\CoInn(H)$. Its kernel consists of those
$\partial(\gamma)$ for which $\gamma=\gamma_1*\gamma_2$ with
$\gamma_1$ a lazy cochain and $\gamma_2$ an algebra morphism
$H\to k$. In other words, the kernel is $B^2_L(H)$.\hfill$\Box$

\smallskip

We shall denote the quotient $\CoInt(H)/\CoInn(H)$ by
$K(H)$. The following Lemma describes $K(H)$ when $H$ is a modified supergroup
algebra.
  
\begin{lemma}Let $H$ be the modified supergroup algebra corresponding
to the group $G$, the representation $V\neq0$, and the central element
$u$. Then $\CoInt(H)\cong{\mathbb Z}_2$, its non-trivial element is
conjugation by $u$ and $K(H)=1$ if and only if
there exists a morphism $\chi\colon G\to k^\cdot$  with $\chi(u)=-1$.
\end{lemma}
\pf Let $\ad(\gamma)\in\CoInt(H)$ with $\gamma(1)=1$. Then $\ad(\gamma)(g)=g$ for
every $g\in G$. Besides, $\ad(\gamma)(gh)=g\,\ad(\gamma)(h)$ and 
$\ad(\gamma)(h)\,g=\ad(\gamma)(hg)$ 
for every
$h\in H$ and every $g\in G$. When $g=u$ this implies that the
${\mathbb Z}_2$-grading induced by conjugation by $u$ must
be preserved by $\ad(\gamma)$. Therefore, since
for every $v\in V$ we have
$\ad(\gamma)(v)=\gamma^{-1}(u)\gamma(v)(u-1)+\gamma^{-1}(u)v$, we
necessarily have 
$\gamma^{-1}(u)\gamma(v)(u-1)=0$ so $\gamma(v)=0$ and  
$\ad(\gamma)(v)=v\gamma^{-1}(u)$ for every $v$. Similarly, 
$\ad(\gamma)(uv)=\gamma(uv)(1-u)+\gamma(u)uv$ must be odd whence
$\gamma(uv)=0$ and 
$\ad(\gamma)(uv)=uv\gamma(u)$. On the other hand,
$\ad(\gamma)(uv)=u\ad(\gamma)(v)$ implies that
$\gamma(u)=\eta$ with
$\eta\in\{\pm1\}$ and necessarily $\ad(\gamma)(gv_{1}\cdots
v_n)=\eta^ngv_{1}\cdots v_n$, that is, $\CoInt(H)\subset\{\id,$
conjugation by $u\}$. 
If we choose representatives $\bar g$ for the
classes of $G/U$ with $1$ representing $U$ we can define
$\gamma(\bar g)=1$, $\gamma(u\bar g)=-1$, and $\gamma(h)=0$ if $h\not\in
k[G]$. Then, conjugation by $u$ is indeed $\ad(\gamma)$ so that
$\CoInt(H)\cong{\mathbb Z}_2$ and we have the first statement. 
Besides, $K(H)$ is trivial if and only if 
conjugation by $u$ is $\ad(\chi)$ for some algebra morphism
$\chi\colon H\to k$. 
This happens if and only if there exists an algebra morphism
$\chi\colon H\to k$ with $\chi(u)=-1$.\hfill$\Box$
%
%

\section{The Brauer group $BM(k, k[G],R_u)$}\label{groupalgebra}

Let $(H, R)$ be a finite dimensional semisimple co-semisimple triangular 
Hopf algebra. 
By \cite{EG3} if $k$ is algebraically closed $(H,R)$ is twist equivalent to
$(k[G], R_u)$ for some finite group $G$ and some
central element $u\in G$ with $u^2=1$. For this reason, we shall
devote this Section to the computation of $BM(k, k[G], R_u)$. 

\smallskip

If $u=1$, then $R_u=1\otimes 1$
so by \cite[Theorem 1.12]{long2} we have: $BM(k, k[G], 1\otimes 1)\cong
Br(k)\times H^2(G, k^\cdot)$. If $u\neq1$  the
role of $H^2(G, k^\cdot)$ will be played by the group which is the
subject of next subsection.

\subsection{The group $H^2_\sharp(G, k^\cdot)$}\label{sharp} 

For a finite group $K$ let $Z^2(K, k^\cdot)$ 
(resp.\ $B^2(K, k^\cdot)$, resp.\ $H^2(K, k^\cdot)$)
denote the group of $2$-cocycles of $K$ (resp.\  the group of $2$
coboundaries of 
$K$, resp.\ the  second cohomology group of $K$ with trivial
coefficients). 

\smallbreak

For any pair of finite groups $L$ and $K$, let
$P(L,K; k^\cdot)$ be the group of pairings of $L$ and $K$ in
$k^\cdot$, i.e., the group of maps $\beta\colon L\times K\to k^\cdot$
such that 
\begin{equation}\label{pairing1}\beta(ab, c)=\beta(a,c)\beta(b,
c)\mbox{ and }\beta(a,cd)=\beta(a,c)\beta(a,d)\end{equation} 
for every $a,b\in L$ and every
$c,d\in K$. It is well-known that $P(L,
K;k^\cdot)\cong\Hom(K,\Hom(L, k^\cdot))\cong \Hom(L,\Hom(K, k^\cdot))$.    

Let $G$, $u\neq1$, $U$ be as in the previous section. The assignment $\theta\colon Z^2(G, k^\cdot)\to P(G, U;
k^\cdot)$ defined by
$\theta(\sigma)(g,u^t)=\sigma(g,u^t)\sigma^{-1}(u^t, g)$ for every
$g\in G$ and $t\in\{0,1\}$, is a 
group morphism inducing a group morphism 
$$\theta\colon H^2(G, k^\cdot)\to{\rm Hom}(G, {\rm Hom}({\mathbb Z}_2, k^\cdot))\cong {\rm
Hom}(G, {\mathbb Z}_2)$$ (see \cite[Lemma
2.2.6]{kar2}).
Hence, each cohomology class $\bar \sigma$ of $G$ determines 
  a ${\mathbb Z}_2$-grading on $k[G]$: the element $g$ is
  even if $\theta(\sigma)(g,u)=1$ and odd if $\theta(\sigma)(g,u)=-1$. 
For all such gradings $u$ is always even. By (\ref{pairing1}) with
$b=u^t$, $c=u^s$
we also see that the image of $\theta$
lies in $P(G/U, U;k^\cdot)\cong\Hom(G/U, {\mathbb Z}_2)$.  

\smallskip

We shall  denote by $|g|_\sigma$ the degree induced by $\sigma$ of the
element $g$, with values in $\{0,1\}$. We point out that for
$\sigma,\omega\in 
Z^2(G, k^\cdot)$ we have
\begin{equation}(-1)^{|g|_{\sigma*\omega}}=\theta(\sigma*\omega)(g,u)=
\theta(\sigma)(g,u) \theta(\omega)(g,u)= (-1)^{
|g|_\sigma+|g|_\omega}.
\end{equation}       
  
Let $Reg^2(G, k^\cdot)$ denote the group of $2$-cochains on $G$.

\begin{proposition}The assignment $\sharp\colon Z^2(G,k^\cdot)\times Z^2(G,k^\cdot)\to
  Reg^2(G, k^\cdot)$  given by $(\sigma\sharp \omega)(g,h):=(\sigma*\omega)(g,h)(-1)^{|g|_\sigma|h|_{\omega}}$ for
  every $\sigma, \omega\in 
  Z^2(G, k^\cdot)$ and every $g,h\in G$, endows $Z^2(G, k^\cdot)$ of a
  group structure inducing an abelian
  group structure on $H^2(G, k^\cdot)$.  
\end{proposition} 
\pf Let $\sigma,\omega\in Z^2(G, k^\cdot)$. Then
$(\sigma\sharp \omega)$ is a $2$-cocycle if and only if  
\begin{equation}\label{gradings}
(-1)^{|g|_\sigma|h|_\omega+|gh|_\sigma|l|_\omega}=(-1)^{|h|_\sigma|l|_\omega+|g|_\sigma|hl|_\omega} 
\end{equation}
for every $g, h,l \in G$. The above equality holds because $|.|_\omega$
and $|.|_\sigma$ are ${\mathbb Z}_2$-gradings, hence $\sharp$ defines
an operation on $Z^2(G,k^\cdot)$.  

Let $\sigma, \omega, \beta\in Z^2(G, k^\cdot)$. Then, for every
$g,h\in G$ we have: 
$$ 
\begin{array}{rl}
((\sigma\sharp \omega)\sharp \beta)(g,h)&=((\sigma*\omega)\sharp \beta)(g,
  h)(-1)^{|g|_\sigma|h|_\omega}\\
&=((\sigma*\omega)*\beta)(g,h)(-1)^{|g|_\sigma|h|_\omega+
|g|_{\sigma*\omega}|h|_\beta}\\ 
&=((\sigma*(\omega*\beta))(g,
  h)(-1)^{|g|_\sigma|h|_\omega+|g|_\sigma|h|_\beta+|g|_\omega|h|_\beta}\\
&=((\sigma*(\omega\sharp \beta)(g, h)(-1)^{|g|_\sigma|h|_\omega+
|g|_\sigma|h|_\beta}\\
&=(\sigma\sharp(\omega\sharp \beta))(g,h)
\end{array}
$$
so $\sharp$ is associative.

The trivial cocycle $\varepsilon\times\varepsilon$  induces the trivial grading on $k[G]$. Hence, 
$$(\sigma\sharp(\varepsilon\times\varepsilon))=(\sigma*(\varepsilon\times\varepsilon))=
\sigma=((\varepsilon\times\varepsilon)*\sigma)=
((\varepsilon\times\varepsilon)\sharp \sigma)$$ and
$\varepsilon\times\varepsilon$ is a neutral element for $(Z^2(G,
k^\cdot),\sharp)$.  

Let $\sigma\in Z^2(G,k^\cdot)$. Then the map $\sigma'\colon G\times G\to
k^\cdot$ defined by
$\sigma'(g,h):=\sigma^{-1}(g,h)(-1)^{|g|_\sigma|h|_\sigma}$ lies in 
$Z^2(G, k^\cdot)$ because (\ref{gradings}) holds when $\sigma=\omega$ and
because $\sigma^{-1}$ is a $2$-cocycle. Since
${\rm Im}(\theta)\subset\Hom(G,{\mathbb Z}_2)$ we have $|g|_\sigma=|g|_{\sigma^{-1}}$so $\sigma\sharp \sigma'=\sigma'\sharp
\sigma=\varepsilon\times\varepsilon$ and $(Z^2(G, 
k^\cdot), \sharp)$ is a group. 

The coboundaries induce the trivial
grading because $u$ is central, hence the $\sharp$ product of an
element in $Z^2(G, k^\cdot)$ with an
element in $B^2(G, k^\cdot)$ coincides
with the usual product. Therefore,
$B^2(G,k^\cdot)$ is a central subgroup of $(Z^2(G,
k^\cdot), \sharp)$ and $H^2(G, k^\cdot)$ inherits the group
structure $\sharp$.

Let now $\sigma,\sigma'\in Z^2(G, k^\cdot)$. Then a priori
 $\sigma\sharp \sigma'\neq 
 \sigma'\sharp \sigma$. However, if we define $\gamma\colon G\to k^\cdot$ by 
$\gamma(g)=(-1)^{|g|_\sigma|g|_{\sigma'}}$ then
 $\gamma(g)\gamma(h)\gamma^{-1}(gh)=(-1)^{|g|_\sigma|h|_{\sigma'}+|g|_{\sigma'}|h|_\sigma}$ 
 so $\sigma\sharp \sigma'$ is cohomologous to $\sigma'\sharp \sigma$
 and $(H^2(G, k^\cdot),\sharp)$ is abelian.\hfill$\Box$

\medskip

From now on we shall denote the group $(H^2(G, k^\cdot), \sharp)$ by
$H^2_\sharp(G, k^\cdot)$. 

\begin{example}{\rm Let $G\cong\Z2\times\Z2$ with generators $x=u$ and $y$. Then $H^2(G, k^\cdot)\cong
  k^\cdot/(k^\cdot)^2\times k^\cdot/(k^\cdot)^2\times\Z2$. Here, the
  first (respectively, second) copy of  $k^\cdot/(k^\cdot)^2$ is represented by cocycles
  $\sigma_a$ (respectively $\omega_a$) such that $\sigma_a(x,x)=a$
  (respectively $\omega_a(y,y)=a$) and trivial elsewhere, for every
  $a\in k^\cdot$. The
  non-trivial class of the copy of $\Z2$ is represented by the cocycle
  $\lambda(x^t y^s, x^ly^m)=(-1)^{tm}$. Then $\sigma_a$ and $\omega_a$
  induce the trivial grading for every $a$ while $\lambda$ induces a
  grading for which $x$ is even and $y$ is odd. Then,
  $(\lambda\sharp\lambda)(x^t y^s, x^ly^m)=(-1)^{sm}$ so
  $\lambda\sharp\lambda=\omega_{-1}$. If $-1$ is a square in $k$ then
  $H^2(G, k^\cdot)\cong H^2_\sharp(G, k^\cdot)$. If $-1$ is not a
  square in $k$ then  $H^2(G, k^\cdot)\not\cong H^2_\sharp(G,
  k^\cdot)$ because there exists an element in  $H^2_\sharp(G,
  k^\cdot)$ whose order is dfferent from $2$. We shall see that this
  example is indeed a key example, in order to see whether the two group
structures on the cohomology  groups are isomorphic.}   
\end{example} 

Next we shall investigate some properties of
$H^2_\sharp(G,k^\cdot)$ that will be needed in the sequel.  Let $\pi\colon G\to G/U$ be the natural
projection. We recall that composition with $\pi\times\pi$ 
induces the inflation morphism 
$\Infl\colon H^2(G/U,k^\cdot)\to H^2(G,k^\cdot)$. The image of $\Infl$
consists of cohomology classes that can be 
represented by doubly $U$-invariant cocycles. 
Composition with $\pi\times\pi$ induces also an inflation morphism 
$H^2(G/U,k^\cdot)\to H_\sharp^2(G,k^\cdot)$, 
which we again
shall denote by $\Infl$, because doubly $U$-invariant cocycles induce
a trivial grading.
The kernel and the image of $\Infl$ coincide respectively 
with the kernel and the image of the classical inflation map.
The classical Hochschild-Serre exact sequence reads:
$$
\begin{array}{l}
\CD
1\longrightarrow\Hom(G/U, k^\cdot)@>\circ\pi>>\Hom(G,k^\cdot)@>{\rm
Res}>>\\
\Hom(U,k^\cdot)@>{\tt T}>>H^2(G/U,k^\cdot)@>{\rm Infl}>>H^2(G, k^\cdot).
\endCD
\end{array}
$$
Here ${\rm Res}$ denotes the usual restriction and
${\tt T}\colon {\rm
  Hom}(U,k^\cdot)\to  H^2(G/U,k^\cdot)$ is the transgression morphism 
defined as follows: if $\phi$ is a section
$G/U\to G$ then  ${\tt
  T}(f)(g,h)=f(\phi(gU)\phi(hU)\phi^{-1}(ghU))$ 
for every $f\in  {\rm
  Hom}(U,k^\cdot)$ and every $g,h\in G$. It follows that ${\rm Infl}$
is injective if and only if ${\tt T}$ is 
trivial if and only if ${\rm Im}({\rm Res})\cong\Hom(U,
k^\cdot)\cong{\mathbb Z}_2$  
if and only if there exists a group morphism $\chi\colon G\to k^\cdot$
such that $\chi(u)=-1$. Thus, the same statement holds for
$H^2_\sharp(G, k^\cdot)$.   

We recall that by \cite[Theorem 2.2.7]{kar2}, there exists
an exact sequence
\begin{equation}\label{read}
\CD
H^2(G/U,k^\cdot)@>{\rm Infl}>>H^2(G, k^\cdot)@>{{\rm
    res}\times\theta}>>H^2(U, k^\cdot)\times\Hom(G/U,{\mathbb Z}_2)
\endCD
\end{equation}
where ${\rm res}$ is the morphism induced by restriction on cocycles. 

If $G\cong U\times G/U$ then ${\rm Infl}$ is injective and ${{\rm
    res}\times\theta}$ is surjective and split, hence 
\begin{equation}\label{SY}
H^2(G, k^\cdot)\cong H^2(U, k^\cdot)\times H^2(G/U,
k^\cdot)\times\Hom(G/U,{\mathbb Z}_2).
\end{equation}
Thus, as a set,
$H^2_\sharp(G, k^\cdot)\cong H^2(U, k^\cdot)\times H^2(G/U,
k^\cdot)\times\Hom(G/U,{\mathbb Z}_2)$. Besides, an analysis of the embedding
of $H^2(G/U,k^\cdot)\times H^2(U, k^\cdot)$ shows that this is a central
subgroup of $H^2_\sharp(G, k^\cdot)$. 
In particular, the product on $H^2_\sharp(G, k^\cdot)$ is given by the
rule:
\begin{equation}\label{rule}
(\overline{\sigma_U},\overline{\sigma_{G/U}},\chi)
\sharp(\overline{\omega_U},\overline{\omega_{G/U}},\chi')=
(\overline{\sigma_U*\omega_U},\overline{\sigma_{G/U}*\omega_{G/U}*
c_{\chi,\chi'}},\chi\chi') 
\end{equation}
where $c_{\chi,\chi'}(gU,hU):=(-1)^{\chi(g)\chi'(h)}$.
 
\smallskip

If $G\not\cong U\times G/U$ sequence (\ref{read}) still yields
an exact sequence for $H^2_\sharp(G,k^\cdot)$:
$$
\begin{array}{l}
\CD
1\longrightarrow\Hom(U, k^\cdot)@>{\tt T}>> H^2(G/U,k^\cdot)@>{\rm
Infl}>>\\
\phantom{1\longrightarrow}H^2_\sharp(G,k^\cdot)@>{{\rm
res}\times\theta}>>
H^2(U,
k^\cdot)\times\Hom(G/U,{\mathbb Z}_2).
\endCD
\end{array}
$$
%
Let us also observe that in this
 case $\Hom(G/U,{\mathbb Z}_2)\cong\Hom(G,{\mathbb Z}_2)$ since there
 is no group morphism $G\to {\mathbb Z}_2\cong U$ taking value
 $-1$ on $u$. 

The following Proposition lists some conditions on $k$ or on $G$
ensuring that $H^2(G,k^\cdot)\cong H^2_\sharp(G, k^\cdot)$. This
allows the reader to compare the results contained in Section \ref{computationG} with those in
\cite{Chi} and \cite{CGO}. By construction, symmetric cocycles induce
the trivial grading so ${\rm Ext}(G, k^\cdot)$ can be seen as a
subgroup of $H^2_\sharp(G, k^\cdot)$. In particular, if $G$ is cyclic,
$H^2_\sharp(G, k^\cdot)=H^2(G, k^\cdot)$. We shall consider some
further cases. 

\begin{proposition}\label{coincide}Let $G, u,$ and $k$ be as before. Then $H^2(G,
 k^\cdot)\cong H^2_\sharp(G,
 k^\cdot)$ in the following cases:
\begin{enumerate}
\item $-1$ is a square in $k$;
\item $G$ is nilpotent, its $2$-Sylow $S_2$ is abelian and  
$u$ is contained in a direct product of cyclic $2$-groups none of which is $\Z2$;
\item $G$ is nilpotent, its $2$-Sylow $S_2$ is abelian and
$S_2\cong\prod_i{\mathbb Z}_{n_i}$ with $n_i=2^{e_i}$ and $e_i\ge2$
  for $i\ge2$.
\end{enumerate} 
\end{proposition}
\pf Let $t(k^\cdot)$
denote the torsion part of
$k^\cdot$. By \cite[Lemma 2.3.19, Theorem 2.3.21]{kar2} we have $H^2(G,
k^\cdot)\cong H^2(G, t(k^\cdot))\times H^2(G, k^\cdot/t(k^\cdot))$ where
all classes in the second summand can be represented by symmetric
cocycles. Hence, $H^2(G, k^\cdot/t(k^\cdot))$  is a subgroup of
$H^2_\sharp(G, k^\cdot)$ and for all its elements $\overline{\sigma}$ and every
element $\overline{\omega}\in H^2_\sharp(G, k^\cdot)$ we have
$\overline{\sigma\sharp\omega}=\overline{\sigma*\omega}$. Moreover,
for every pair of cocycles $\omega$ and $\sigma$ in $Z^2(G,
t(k^\cdot))$, their product $\sigma\sharp\omega$ and their inverses
still take values in $t(k^\cdot)$, so  $H^2_\sharp(G,
t(k^\cdot))$ is a subgroup of $H^2_\sharp(G, k^\cdot)$ and
$H^2_\sharp(G, k^\cdot)\cong H^2_\sharp(G,
t(k^\cdot))\times H^2(G, k^\cdot/t(k^\cdot))$.

Let us assume that  $-1=\zeta^2$ for some $\zeta\in k$. For every morphism
$\chi\colon G/U\to \Z2$ we define $\mu\colon k[G]\to
k^\cdot$ by $$\mu(g)=\cases{1& if $\chi(g)=1$\cr
\zeta& if $\chi(g)=-1$\cr}$$ 
and we put $|g|_\chi=0$ if $\chi(g)=1$ and $|g|_\chi=1$ if $\chi(g)=-1$. 
Then we have: 
$$
\begin{array}{rl}
\partial(\mu)(g,h)&=\mu(g)\mu(h)\mu^{-1}(gh)\\
&=\zeta^{|g|_\chi}\zeta^{|h|_\chi}(\zeta^{-1})^{|gh|_\chi}\\
&=\cases{1& if $\chi(g)=1$ or $\chi(h)=1$\cr
\zeta^2 & if $\chi(g)=\chi(h)=-1$\cr}\\
&=(-1)^{\chi(g)\chi(h)}.
\end{array}
$$
Thus, for every cocycle $\sigma$, the products $\sigma\sharp\sigma$
and $\sigma*\sigma$ are cohomologous. It follows that, for every class
$\overline\sigma$, the order of $\overline\sigma$ is the same in
$H^2_\sharp(G, k^\cdot)$ and in $H^2(G, k^\cdot)$. Since $H^2_\sharp(G,
t(k^\cdot))$ is a finite abelian group of the same order as $H^2(G,
t(k^\cdot))$ and since for every $m$ the number of elements of order
$m$ in both groups coincide, we have the statement in case $1$. 

\smallbreak

Let us now assume that $G\cong \prod_{p}S_p$ is the direct product of
its Sylow subgroups. By \cite[Corollary 2.3.15]{kar2}, 
$H^2(G,k^\cdot)\cong \prod_{p}H^2(S_p,k^\cdot)$. 
The cocycle $\sigma$ representing an element in
$H^2(S_p, k^\cdot)$ with $p$ odd induces a trivial
grading so $\sigma\sharp \omega=\sigma*\omega$ for every $\omega\in
Z^2(G, k^\cdot)$. 
Besides, the $\sharp$-product of two cocycles representing
elements in $H^2(S_2, k^\cdot)$ is represented by a 2-cocycle on
$S_2\times S_2$ because $\theta(\sigma)(h,u)=1$ for every 
$h\in\prod_{p\neq2} S_p$. Therefore, we have a group isomorphism 
$$H^2_\sharp(G, k^\cdot)\cong
H^2_\sharp(S_2,k^\cdot)\times\prod_{p\neq2}H^2(S_p,k^\cdot)$$ 
%
and it is enough to show that $H^2(S_2, k^\cdot)\cong
H^2_\sharp(S_2,k^\cdot)$. Let us assume that we are not in case $1$
and that $S_2$ is abelian, so that
$S_2\cong\prod_j{\mathbb Z}_{n_j}$ with 
$n_j=2^{e_j}$ and $e_i\le e_{i+1}$. Then \cite[Theorem 2.3.13]{kar2} implies that 
$$H^2(S_2,k^\cdot)\cong
\prod_jk^\cdot/(k^\cdot)^{n_j}\times\prod_{i<j}P({\mathbb
  Z}_{n_i},{\mathbb Z}_{n_j}; k^\cdot).$$ The elements in
$k^\cdot/(k^\cdot)^{n_j}$ are represented by symmetric cocycles for
every $j$, hence these groups are also subgroups of
$H^2_\sharp(S_2,k^\cdot)$. Since 
$-1$ is not a square $P({\mathbb
  Z}_{n_i},{\mathbb Z}_{n_j}; k^\cdot)\cong \Hom({\mathbb
  Z}_{n_i},\Hom({\mathbb Z}_{n_j}, {\mathbb Z}_2))\cong \Hom({\mathbb
  Z}_{n_i}, {\mathbb Z}_2)\cong \Z2$. If we put ${\mathbb
  Z}_{n_i}\cong\langle v_i\rangle$ for every $i$, the non-trivial element in 
$P({\mathbb
  Z}_{n_i},{\mathbb Z}_{n_j}; k^\cdot)$ is represented by
  $\beta_{ij}(v_i^s,v_j^t)=(-1)^{st}$. The corresponding
  cohomology class is represented by the cocycle $\sigma_{ij}$ on ${\mathbb
  Z}_{n_i}\times{\mathbb Z}_{n_j}$ such that
  $\sigma_{ij}(v_i^ev_j^f,v_i^lv_j^t)=\beta_{ij}(v_i^e,v_j^t)$ and
  trivial on the remaining summands. Let us
  determine $\theta(\sigma_{ij})$. The element $u$ can be written as
 $\prod_t v_t^{l_t}$ with $l_t\in\{0, 2^{e_t-1}\}$. For every
  $g=v_1^{s_1}\cdots v_m^{s_m}$ we have
$$\sigma_{ij}(v_1^{s_1}\cdots v_m^{s_m},
  u)\sigma^{-1}_{ij}(u,v_1^{s_1}\cdots
  v_m^{s_m})=(-1)^{s_il_j+l_is_j}.$$
If $u$ is as in case $2$ or if no summand isomorphic to $\Z2$ occurs
  in $S_2$ then $l_t$ is even for
  every $t$, the grading is always trivial and we have the statement.  

Let us then assume that $e_1=1$, $l_1=1$ and that there is one
 summand isomorphic to ${\mathbb Z}_2$ in $S_2$. Then
 $\sigma_{jk}\sharp \sigma_{rs}=\sigma_{jk}*\sigma_{rs}$ if $r,s,j,k\neq1$ while  
 we have $$\theta(\sigma_{1j})(v_1^{s_1}\cdots
  v_m^{s_m},u)=(-1)^{s_j}.$$ Thus, for $j\le t$ we have
 $(\sigma_{1j}\sharp
  \sigma_{1t})=(\sigma_{1j}*\sigma_{1t})*c_{jt}$ where 
$$c_{jt}(v_1^{s_1}\cdots
  v_m^{s_m},v_1^{p_1}\cdots
  v_m^{p_m})=(-1)^{s_jp_t}.$$ 
When $j\neq t$, the cocycle 
$c_{jt}$ is $\sigma_{jt}$ introduced before. 
If $j=t$ and $n_j=2^{e_j}$ with $e_j\ge2$ then 
  $c_{jj}$ is a coboundary. Indeed, the cochain $\mu$ given by
 $$\mu(\prod_{b=1}^mv_b^{a_b})=(-1)^{{a_j}\choose2}$$ is well-defined on $S_2$ and
$c_{jj}=\partial(\mu)$. It follows that the set $\prod_{i<j}P({\mathbb Z}_{n_i},{\mathbb Z}_{n_j};k^\cdot)$ forms a
  subgroup of $H^2_\sharp(S_2, k^\cdot)$ whose elements have all order
  $\le 2$. Hence, it must be isomorphic to $\prod_{i<j}P({\mathbb
    Z}_{n_i},{\mathbb Z}_{n_j};k^\cdot)$ with the usual product. Thus,  $H^2_\sharp(S_2,
  k^\cdot)\cong \prod_j k^\cdot/(k^\cdot)^{n_j}\times
  \prod_{i<j}{\mathbb Z}_2$ and we have the statement in case $3$. 
\hfill$\Box$
%

\begin{corollary}\label{closed}If $k=\bar k$ then $H^2_\sharp(G, k^\cdot)\cong H^2(G, k^\cdot)$.\hfill$\Box$ 
\end{corollary}

\subsection{The computation of $BM(k, k[G], R_u)$} \label{computationG}

In this Section we shall compute $BM(k, k[G], R_u)$ when
$u\neq1$.

The pull-back along the embedding $i\colon{k}[U]\to k[G]$ induces a
group morphism $i^*\colon BM(k, k[G], R_u)\to BM(k, k[U], R_u)\cong BW(k)$. As a set, $BW(k)\cong
Br(k)\times k^\cdot/(k^\cdot)^2\times{\mathbb Z}_2$. The elements of
$Br(k)$ are represented by central simple algebras with trivial
${\mathbb Z}_2$-grading. The non-trivial element of ${\mathbb Z}_2$ is
represented by the algebra $C(1)$ generated by the odd element 
$x$ with relation $x^2=1$. The elements of $k^\cdot/(k^\cdot)^2$ are
represented by algebras $C(1)\sharp C(-\alpha)$  with $\alpha\in
k^\cdot/(k^\cdot)^2$, where $C(-\alpha)$ is  generated by the odd element 
$y$ with the relation $y^2=-\alpha$ and with $\sharp$ as in Section
\ref{definizioni}. More precisely, 
to the class of a ${\mathbb Z}_2$-graded central
simple algebra $A$ one associates: $([A],1,0)$ if its odd part
$A_1=0$; $([A],f(u)^2,0)$ if $A$ is central simple and $f\colon
k[U]\to A$ is the map realizing the ${\mathbb Z}_2$-action;
$([A_0],\delta(A),1)$ if 
$A$ is not central simple, $A_0$ is its even part, $\delta(A)=z^2$ where $z\in Z(A)$,
$z\not\in k$ and $z^2\in k$ (\cite[Theorem V.3.10]{lam}). Let us analyze
the image of $i^*$.   

\begin{lemma}\label{surjective}The following assertions are equivalent:
\begin{enumerate}
\item The morphism $i^*$ is surjective;
\item The subgroup $U$ is a direct summand of $G$;
\item The morphism $i^*$ is surjective and split.
\end{enumerate}
\end{lemma}
\pf If $i^*$ is surjective then $[C(1)]$ lies in the image of $i^*$,
that is, there exists a $k[G]$-module algebra $A$ which is ${\mathbb
Z}_2$-graded central simple and isomorphic, as
a $k[U]$-module algebra, to $C(1)\sharp {\rm \End}(P)$ for some
$k[U]$-module $P$, with $U$-action on ${\rm End}(P)$ induced by the
$U$-action on $P$ as in (\ref{endo}). Since the $U$-action on ${\rm End}(P)$ is strongly
inner, $C(1)\sharp {\rm End}(P)\cong C(1)\otimes{\rm End}(P)$ as a
$k[U]$-module algebra. An isomorphism is given by $a\sharp F\mapsto
a\otimes f(u)^{|a|}F$ if $f$ is the algebra morphism $k[U]\to{\rm End}(P)$
realizing the $U$-action (\cite[Theorem IV.2.7]{lam}). The centre of $A$ is a $k[G]$-submodule
algebra. 
Indeed, if $z\in Z(A)$ then for every $y\in A$ and every $g\in G$ 
$$(g.z)y=(g.z)(g.g^{-1}.y)=g.(z (g^{-1}.y))=g.((g^{-1}.y) z)=y(g.z).$$
Besides, $Z(A)$ is isomorphic to $C(1)$ as a $k[U]$-module
algebra because $A\cong C(1)\otimes\End(P)$ as an algebra. Thus, if $i^*$ is surjective then we can lift the 
$U$-action on $C(1)$ to a $G$-action.  If this is
possible, for every $g\in G$ we have $g.x=a_gx+b_g$ for some scalars
$a_g, b_g$ and $ug.x=gu.x=-g.x$
implies that $b_g=0$ for every $g$. 
In this case, $g.x=\chi(g)x$ for some algebra
morphism $\chi\colon{k}[G]\to{k}$ with $\chi(u)=-1$.
Besides, since $x^2=1$, the $k[G]$-module algebra condition for
$g.x^2$ implies that
$\chi(g)^2=1$. Hence if $i^*$ is surjective there exists a group
morphism $\chi\colon G\to {\mathbb Z}_2\cong\{\pm 1\}$ with
$\chi(u)=-1$. Then $G\cong U\times G/U$ so 1 implies 2. 

If 2 holds, the pull-back along the projection map $G\to U$ induces a group
morphism splitting $i^*$ at the Brauer group level, so 2 implies 1 and
3.    
\hfill$\Box$
  
\bigskip

\begin{remark}\label{alfa}{\rm Let us observe that the proof of Lemma
    \ref{surjective} 
shows also that $[C(\alpha)]\in{\rm
Im}(i^*)$ if and only if $U$ is a direct summand of $G$. }
\end{remark} 

Let us denote by ${\rm res}\colon H^2(G, k^\cdot)\to H^2(U, k^\cdot)\cong
k^\cdot/(k^\cdot)^2$ 
the morphism induced by the restriction map. We recall that $BW(k)$ is
a central extension of $Br(k)$ by a group $Q(k)$ which is
$k^\cdot/(k^\cdot)^2\times\Z2$ as a set, with multiplication rule
$(\bar{\alpha},(-1)^e)(\bar{\beta}, (-1)^f)=(\overline{\alpha\beta}(-
1)^{ef},(-1)^{e+f})$ with $e,\, f\in\{\pm 1\}$. In particular, the
central extension of $Br(k)$ by ${\rm res}(H^2(G, k^\cdot))$ is a
subgroup of $BW(k)$, which we shall denote by $Br(k)*{\rm res}(H^2(G,
k^\cdot))$.      

\begin{proposition}\label{imma}If $U$ is not a direct summand of $G$
  the image of 
the map $i^*$ is $Br(k)*{\rm res}(H^2(G, k^\cdot))$.
\end{proposition}
\pf Any central simple algebra $B$ with trivial
$G$-action represents an element $[B]$ in $BM(k, k[G], R_u)$ whose image
under $i^*$ is exactly $[B]\in Br(k)$. Hence $Br(k)$ lies in the image
of $i^*$. It follows from the proof of Lemma \ref{surjective} and Remark
\ref{alfa}  that if
$U$ is not a direct summand of $G$, the classes $[C(1)]$ and
$[C(\alpha)]$ do not lie 
in ${\rm Im}(i^*)$. The proof of the proposition follows from the knowledge of the
multiplication rules in $BW(k)$ (\cite[Theorem V.3.9]{lam}) if we show that
$[C(1)\sharp C(-\alpha)]\in{\rm Im}(i^*)$ if and only if $\overline{\alpha}\in{\rm
Im}({\rm res})\subset k^\cdot/(k^\cdot)^2$. 

Let $[C(1)\sharp C(-\alpha)]\in{\rm Im}(i^*)$. Then for some
$k[G]$-module algebra $B$ which is ${\mathbb Z}_2$-graded central
simple, $B\cong C(1)\sharp C(-\alpha)\sharp{\rm End}(P)$ as
$k[U]$-module algebras. As above, $B\cong (C(1)\sharp
C(-\alpha))\otimes{\rm End}(P)$, so $B$ is central simple (and in fact,
a matrix algebra). Then the
$k[G]$-action on $B$ is necessarily inner and if $f\colon k[G]\to B$
is the map realizing the action, $f(g)f(h)=f(gh)c(g,h)$ for some
2-cocycle $c$ for $G$. The restriction of $c$ to $U\times U$ is
determined by $c(u,u)=\beta$. It is not hard to see
that the (cohomology) class of ${\rm res}(c)$ does not depend on the
choice of the representative of $[C(1)\sharp C(-\alpha)]\in BW(k)$. Let us
choose then $C(1)\sharp C(-\alpha)$. The element realizing, by
conjugation, the $u$-action on $C(1)\sharp C(-\alpha)$  is $x\sharp
y$, where $x^2=1$ and 
$y^2=-\alpha$ are odd generators of $C(1)$ and $C(-\alpha)$,
respectively. Then $(x\sharp y)^2=\alpha$ is in the same class modulo
$(k^\cdot)^2$ as $\beta$. Therefore, if $[C(1)\sharp C(-\alpha)]\in{\rm
Im}(i^*)$ then  $\bar{\alpha}\equiv{\rm res}(c)(u,u)\;{\rm
  mod}(k^\cdot)^2$ for some 
 $c\in Z^2(G, k^\cdot)$.

Conversely, let $\bar\alpha\in k^\cdot/(k^\cdot)^2$ lie in the image of
${\rm res}$. Then there exists a 2-cocycle $c$ such that
$c(u,u)=\alpha$. We may construct as in \cite[Page 567]{long2} 
the $k[G]$-module algebra $A^c$ with underlying algebra ${\rm
  End}(k[G])$. We recall that the $G$-action on $F\in A^c$ is given by
$g.F=c^{-1}(g^{-1}, g)f_g\circ
F\circ f_{g^{-1}}$ where $f_g(h)=c(g,h)gh$. 
If we prove that $A^c$ is $R_u$-Azumaya, then
$i^*([A^c])$ corresponds to the triple $([k],\overline{c(u,u)},0)$
i.e., $i^*([A^c])=[C(1)\sharp C(-\alpha)]$.  

The algebra $A^c$ will be 
$R_u$-Azumaya if and only if it is ${\mathbb Z}_2$-graded central
simple. This happens if and only if its ${\mathbb Z}_2$-graded centre is
trivial because $A^c$ is simple. The $u$-action on $A^c$ is given
by:$(u.F):=\frac{1}{c(u,u)}f_u\circ F\circ f_u$.  

Let $\mu\colon G/U\to G$ be a section for $G\to G/U$, with $\mu(U)=u$. 
A basis for the even part $A^c_0$ of $A^c$ is given by the elements
$$F_{g,h}:=h\otimes g^*+\frac{c(u,h)}{c(u,g)}hu\otimes (gu)^*$$ for $g\in G$
and  $h\in\mu(G/U)$. A basis for the odd part $A^c_1$ of $A^c$ is
given by the elements 
$$L_{g,h}:=h\otimes g^*-\frac{c(u,h)}{c(u,g)}hu\otimes (gu)^*$$ for $g\in G$
and  $h\in\mu(G/U)$. 

Let $z$ be an element of the graded centre, with 
homogeneous summands $z_0=\sum_{g,h}a_{gh}F_{gh}$ and $z_1=
\sum b_{gh}L_{gh}$. 

Commutation of $z_0$ with $F_{r,s}$ evaluated in $r\in \mu(G/U)$ shows that
$a_{gh}=0$ unless $h\in\{g,gu\}$ and that
$$a_{rr}=a_{ss}\quad\mbox{\rm and}\quad
\frac{a_{su,s}}{c(u,su)}=\frac{a_{ru,r}}{c(u,ru)}\quad\forall r,s\in\mu(G/U)$$
so $z_0=a_{u,u}{\rm
  Id}+a_{1,u}\sum_{g\in\mu(G/U)}c(u,gu)F_{gu,g}$. For $r,s\in\mu(G/U)$
we have
$0=z_0L_{rs}-L_{rs}z_0=2a_{1,u}c(u,s)su$, hence $z_0$ is a multiple of
the identity. Commutation of $z_1$ with $F_{r,s}$
shows that $b_{gh}=0$ unless $g\in\{h, uh\}$, that $b_{ss}=b_{rr}$ and
that $b_{su,s}=-b_{ru,r}\frac{c(u,su)}{c(u,ru)}$. Commutation of $z_1$
with $F_{ru,s}$ shows that $b_{ss}=-b_{rr}$ and that $b_{su,s}=b_{ru,r}\frac{c(u,su)}{c(u,ru)}$.
Thus, $z_1=0$ and we have the statement.\hfill$\Box$
%

\smallskip

\begin{proposition}\label{subgroup}If $G\not\cong G/U\times U$
  each element in 
 $BM(k,k[G], R_u)$ can be represented by a central simple algebra.
 If $G\cong U\times G/U$ each element in
  $BM(k,k[G], R_u)$ can be represented by a 
product of the form
  $B\sharp C(1)^a$ with $B$ central simple and $a=0,1$. In particular, 
 the elements that can be represented by central simple algebras form 
 a subgroup of $BM(k, k[G], R_u)$.  
\end{proposition} 
\pf Let $[M]$ be an element in $BM(k,k[G],R_u)$. The braiding induced by
$R_u$ shows that an algebra is $(k[G],R_u)$-Azumaya if and only if it
is $\Z2$-graded central simple.
If $U$ is not a direct summand of $G$ then $i^*([M])=[A\sharp
  C(1)\sharp C(-\alpha)]$ in $BW(k)$ for some central simple algebra
$A$ with trivial action and for some scalar $\alpha$. Thus,  for some
$k[U]$-module $P$ we have
$k[U]$-module algebra isomorphisms
$$
\begin{array}{rl}
M&\cong \End(P)\sharp A\sharp C(1)\sharp C(-\alpha)\\
&\cong\End(P) \otimes (A\sharp C(1)\sharp C(-\alpha))\\
&\cong\End(P)\otimes A\otimes (C(1)\sharp C(-\alpha))
\end{array}
$$ 
so $M$ is central simple. 

If $U$ is a direct summand of $G$ then $i^*([M])=[A\sharp 
  C(1)\sharp C(-\alpha)\sharp C(1)^a]$ in $BW(k)$  for some central simple algebra
$A$ with trivial action, some scalar $\alpha$ and some $a\in\{0,1\}$. 
Thus,  for some
$k[U]$-module $P$ we have
$k[U]$-module algebra isomorphisms
$$
\begin{array}{rl}
M&\cong \End(P)\sharp A\sharp C(1)\sharp C(-\alpha)\sharp C(1)^a\\
&\cong (\End(P)\otimes A\otimes (C(1)\sharp C(-\alpha)))\sharp C(1)^a\\
\end{array}
$$ 
so $M$ is of the required form. The elements represented by a central
simple algebra are those with $a=0$. They form a subgroup, namely the
preimage of the subgroup of $BW(k)$ consisting of triples
$([A],\alpha, 0)$ (for the multiplication rules in $BW(k)$, see \cite[Theorem V.3.9(2)]{lam}).\hfill$\Box$
\bigskip  

We are ready to state the main results of this Section.

\begin{theorem}\label{notdirect}If $G\not\cong G/U\times U$ there is a
  short exact sequence 
$$
1\longrightarrow Br(k)\longrightarrow BM(k, k[G], R_u)\longrightarrow
H^2_\sharp(G, k^\cdot)\longrightarrow 1.$$
If $-1$ is a square in $k$ then $H^2_\sharp(G, k^\cdot)$ can be
replaced by $H^2(G, k^\cdot)$. 
\end{theorem}
\pf By Proposition \ref{subgroup} each element in $BM(k,k[G], R_u)$
can be represented by a central 
simple algebra $B$. The action of $G$ is inner and if
$f\colon k[G]\to B$ is a convolution invertible map for which
$g.b=f(g)bf(g)^{-1}$ then any other such map is of the form $f*\gamma$
with $\gamma$ a convolution invertible map $k[G]\to k$. Therefore the
subalgebra generated by $f(g)$ for all $g\in G$ is well-defined and 
since $B$ is central, $f(g)f(h)=f(gh)c(g,h)$ for some $c\in
Z^2(G, k^\cdot)$. A different choice of the map $f$ would yield a
cohomologous cocycle, so we may associate a cohomology class to $B$. 
If we consider a different representative for $[B]$, say, $B\sharp
{\rm End}(P)$ for some $k[G]$-module $P$, then we have an algebra map
$f'\colon k[G]\to {\rm End}(P)$ realizing the $k[G]$-action on ${\rm
  End}(P)$ and   
$$g.(b\sharp F)=g.b\sharp g.F=f(g)bf^{-1}(g)\sharp f'(g)\circ F\circ
(f'(g))^{-1}.$$
Since $f'$ is an algebra map, each $f'(g)$ is even with respect to
the ${\mathbb Z}_2$-grading induced on ${\rm End}(P)$ by the action of
$u$. Then for every $g,h\in G$ the inverse of $f(h)\sharp
f'(g)$ in $B\sharp {\rm End}(P)$ is $f(h)^{-1}\sharp f'(g)^{-1}$. Besides, since
$$f(g)f(u)=f(gu)c(g,u)=f(ug)c(u,g)\theta(c)(g,u)=f(u)f(g)\theta(c)(g,u)$$
the parity of $g$ with respect to the ${\mathbb Z}_2$-grading
induced by $c$ coincides with the parity of $f(g)$ with
respect to the ${\mathbb Z}_2$-grading induced bu the $u$-action on
$B$ and the same holds for $f'$, $\End(P)$. 

Since the $u$-action is strongly
inner, $B\sharp {\rm End}(P)\cong B\otimes{\rm End}(P)$ as an algebra, 
hence it is central simple and the action is again inner. 

Then, for homogeneous $b\in B$ and $F\in {\rm End}(P)$ we have:
$$
\begin{array}{l}
(f(g)\sharp f'(gu^{|g|_c}))(b\sharp F)(f^{-1}(g)\sharp
  f'(gu^{|g|_c})^{-1})\\
=(f(g)b\sharp f'(gu^{|g|_c})F)(f^{-1}(g)\sharp
  f'(gu^{|g|_c})^{-1}\\
=(-1)^{|F||g|_c}f(g)bf^{-1}(g)\sharp f'(gu^{|g|_c})\circ F\circ
  f'(gu^{|g|_c})^{-1}\\    
=(-1)^{|F||g|_c}g.b\sharp (-1)^{|g|_c|F|}g.F= g.(b\sharp F)
\end{array}
$$
so the map $\bar f\colon k[G]\to B\sharp {\rm End}(P)$ given by $g\mapsto f(g)\sharp f'(gu^{|g|_c})$
realizes the $G$-action.  
Then 
$$
\begin{array}{rl}
\bar f(g)\bar f(h)&=((f(g)\sharp f'(gu^{|g|_c}))(f(h)\sharp f'(hu^{|h|_c}))\\
&=f(gh)c(g,h)\sharp f'(gu^{|g|_c})f'(hu^{|h|_c})\\
&=f(gh)\sharp f'(ghu^{|g|_c+|h|_c})c(g,h)\\
&=\bar f(gh)c(g,h).
\end{array}
$$ 
Thus, the cohomology class of $c$ does not depend on the choice of the 
representative of $[B]$ and we have a well-defined map $BM(k, k[G],
R_u)\to H^2(G, k^\cdot)$.  

Let $[B]$ and $[C]\in BM(k, k[G], R_u)$ with maps $f$ and $s$ inducing
the 
action and with corresponding cocycles $d$ and $e$. We wish to relate
the cohomology class corresponding to $B\sharp C$ with $\bar d$ and
$\bar e$.

The degree of $f(g)$ (resp. $s(g)$) in $B$
(resp. $C$) for the grading induced by the $u$-action coincides with
the degree of
$g$ induced by the cocycle $d$ (resp. $e$) through the map $\theta$. 
It is not hard to verify that 
$$(f(g)f(u)^{|g|_e}\sharp s(g)s(u)^{|g|_d})^{-1}=
(-1)^{|g|_e|g|_d}f^{-1}(u)^{|g|_e}f^{-1}(g)\sharp s^{-1}(u)^{|g|_d}s^{-1}(g).
$$
Let $b\in B$, $c\in C$ be homogeneous. Then  
$$
\begin{array}{l}
(f(g)f(u)^{|g|_e}\sharp s(g)s(u)^{|g|_d})\times\\
\phantom{(f(g)f(u)}(b\sharp
c)((-1)^{|g|_e|g|_d}f^{-1}(u)^{|g|_e}f^{-1}(g)
\sharp s^{-1}(u)^{|g|_d} s^{-1}(g))\\
=(-1)^{|b|_d|g|_e+|g|_d|g|_e}
(f(g)f(u)^{|g|_e}b\sharp s(g)s(u)^{|g|_d}c)\times\\
\phantom{(f(g)f(u)}(f^{-1}(u)^{|g|_e}f^{-1}(g)\sharp
 s^{-1}(u)^{|g|_d} s^{-1}(g))\\
=(-1)^{|b|_d|g|_e+|c|_e|g|_d}\\
\phantom{(f(g)f(u)}
(f(g)f(u)^{|g|_e}bf^{-1}(u)^{|g|_e}f^{-1}(g)\sharp s(g)s(u)^{|g|_d}c
 s^{-1}(u)^{|g|_d} s^{-1}(g))\\
=(f(g)bf^{-1}(g)\sharp s(g)c
 s^{-1}(g))\\
=g.b\sharp g.c
\end{array}
$$
hence the map $\bar f\colon k[G]\to B\sharp C$ given by $\bar
f(g)=(f(g)f(u)^{|g|_e}\sharp s(g)s(u)^{|g|_d})$ realizes the
$G$-action on $B\sharp C$. 
The computation
$$
\begin{array}{rl}
\bar f(g)\bar f(h)&=(f(g)f(u)^{|g|_e}\sharp
s(g)s(u)^{|g|_d})(f(h)f(u)^{|h|_e}\sharp
s(h)s(u)^{|h|_d}\\
&=(-1)^{2|h|_d|g|_e+|h|_e|g|_d}
(f(g)f(h)f(u)^{|gh|_e}\sharp
s(g)s(h)s(u)^{|gh|_d})\\
&=(-1)^{|g|_d|h|_e}f(gh)f(u)^{|gh|_e}d(g,h)\sharp
s(gh)s(u)^{|gh|_d}e(g,h)\\
&=(d\sharp e)(g,h)\bar f(gh)
\end{array}
$$
shows that the assignment $\zeta\colon BM(k,k[G],R_u)\to H^2_\sharp(G,
k^\cdot)$ defined by $[B]\mapsto\bar{d}$ 
is a group morphism. The construction of
$A^c$ in the proof of Proposition \ref{imma} shows that $\zeta$ is
surjective. 
The kernel of the morphism $\zeta$ contains $Br(k)$ because every
central simple algebra with trivial $G$-action represents an element
in the kernel of $\zeta$. Let $[B]\in{\rm Ker(\zeta)}$. Then $B$ is central simple, with strongly inner $G$-action. 
Let us consider the algebra $B^{op}_{triv}$, which is isomorphic to $B^{op}$ as an
algebra, and with trivial $G$-action. Then, as algebras, 
$$B\sharp B^{op}_{triv}\cong B\otimes B^{op}_{triv}\cong B\otimes B^{op}\cong{\rm End}(B).$$  
Besides, if $f$ is the algebra map realizing the $G$-action on $B$,
then $g\lu(b\sharp b')=(f(g)\sharp 1)(b\sharp b')(f(g)^{-1}\sharp 1)$ so
$B\sharp B^{op}_{triv}$ is a matrix algebra with strongly inner
$k[G]$-action. Thus, $[B]=[B^{op}_{triv}]^{-1}=[B_{triv}]$ 
can be represented by a central simple 
algebra with trivial action and ${\rm Ker}(\zeta)=Br(k)$. The last
statement follows from Corollary \ref{closed}.\hfill$\Box$  

\smallskip

\begin{theorem}\label{sidirect}If $G\cong U\times G/U$ there is a short exact sequence
$$
1\longrightarrow Br(k)\longrightarrow BM(k, k[G], R_u)\longrightarrow
Q(k,G)\longrightarrow 1$$
where $Q(k,G)$ is the group with underlying set
$$H^2_\sharp(G,
    k^\cdot)\times{\mathbb Z}_2\cong H^2(G/U, k^\cdot)\times
    \Hom(G/U,{\mathbb Z}_2)\times k^\cdot/(k^\cdot)^2\times {\mathbb Z}_2$$
and with multiplication rule
$$
\begin{array}{l}
(\overline{\sigma_{G/U}},\chi,\bar{t},(-1)^e)(\overline{\omega_{G/U}},\chi',\bar{s},(-1)^f)\\
\phantom{(\bar\sigma_{G/U},}=(\overline{\sigma_{G/U}*\omega_{G/U}*c_{\chi
    \chi'}}, \chi\chi', \overline{st}(-1)^{ef}, (-1)^{e+f})
\end{array}
$$ with $c_{\chi, \chi'}(g,h)=(-1)^{\chi(g)\chi'(h)}$. 

If $-1$ is a square in $k$ we may replace $Q(k,G)$ by $H^2(G,
    k^\cdot)\times{\mathbb Z}_2$. 
\end{theorem}
\pf By Proposition \ref{subgroup} each class of $BM(k,k[G],R_u)$ is
represented by some product $B\sharp C(1)^a$ with $B$ central simple
and $a\in\{0,1\}$. As in the
previous case we may associate a cohomology class of $G$ to each class
represented by a central simple algebra.
%
%
Let us define the map
$$
\begin{array}{rl}
\zeta\colon BM(k, k[G],
R_u)&\longrightarrow Q(k, G)\\
\left[ B\sharp C(1)^a\right]& \mapsto (\overline{\sigma}_B, (-1)^a)
\end{array}
$$
where $\overline{\sigma}_B$ is the cohomology class associated with $B$.
This is a well-defined map because 
the cohomology class, centrality and simplicity of an algebra do not
depend 
on the choice of the representative.

The map $\zeta$ is a group morphism. Indeed, if two classes and  are
represented by a central simple algebra then the image lies in
$H^2_\sharp(G,k^\cdot)\times 1$ and the proof is as in Theorem \ref{notdirect}.
If we have $B$ and $C\sharp C(1)$, with $B$ and $C$ central simple then
 $B\sharp C\sharp C(1)$ is not central simple. Then
$$\zeta([B\sharp C\sharp C(1)])=(\overline{\sigma}_{B\sharp C}, -1)=
(\overline{\sigma}_B\sharp\overline{\sigma}_C,-1)=
\zeta([B])\zeta([C\sharp C(1)])$$
follows as in Theorem \ref{notdirect}. 
If we are given $B\sharp C(1)$ and $C\sharp C(1)$ with $B$ and $C$
central simple, then $A:=B\sharp C(1)\sharp C\sharp C(1)\cong B\sharp
C\sharp C(1)\sharp C(1)$ is central
simple by \cite[Theorem V.3.9]{lam}. 
If $f$ is the map
realizing the action on $B\sharp C$ then $f(u)\sharp(x\sharp y)$
with $x, y$ homogeneous generators of $C(1)$
realizes the $u$-action on $A$, while $f(g)\sharp 1$ realizes the
$g$-action on $A$ for every $g\in G/U$.
Therefore if
$\sigma_{B\sharp C}=(\sigma_{G/U},\chi,\alpha)$ then
$\sigma_A=(\sigma_{G/U}, \chi,-\alpha)$. Hence $\zeta([B\sharp
  C(1)][C\sharp C(1)])=\zeta([B\sharp C(1)])\zeta([C\sharp
  C(1)])$ and $\zeta$ is a group morphism. 
As is Theorem \ref{notdirect}, it is not hard to show that ${\rm
 Ker}(\zeta)\cong Br(k)$. Let $(\overline{\sigma},(-1)^e)\in Q(k, G)$. 
Then $A^\sigma\sharp C(1)^e$ represents an
 element in $BM(k,k[G], R_u)$ whose image under $\zeta$ is
 $(\overline{\sigma},(-1)^e)$. Hence $\zeta$ is bijective and we have the
first statement. The second statement follows from Corollary \ref{closed}.\hfill$\Box$  

\begin{remark}{\rm If $G\cong U\times G/U$ then
 $BM(k,k[G],R_u)$ can also be given as a central extension of $Br(k)$ with
a group isomorphic to $H^2(G/U,k^\cdot)\times
H^2(U,k^\cdot)\times\Hom(G,{\mathbb Z}_2)$ as a set and with
multiplication
 rule
\begin{equation}\label{rule1}
(\overline{\sigma_U},\overline{\sigma_{G/U}},\chi)
\sharp(\overline{\omega_U},\overline{\omega_{G/U}},\chi')=
(\overline{\sigma_U*\omega_U},\overline{\sigma_{G/U}*\omega_{G/U}*c_{\chi,\chi'}},\chi\chi') 
\end{equation}
where $c_{\chi,\chi'}(g,h):=(-1)^{\chi(g)\chi'(h)}$.}
\end{remark}

In Section \ref{comparison} we shall compare the result in this and
the following sections with classical exact sequences involving Brauer
goups and groups of Galois objects. 

\medskip

It follows from \cite[Theorem 5.1, Theorem 5.2, Corollary 6.2]{EG3}
that all semi-simple and cosemi-simple triangular Hopf algebras over
an algebraically closed field are Drinfeld twists of group algebras
with $R=1\otimes 1$ or $R=R_u$ for some central element $u$ of order
$2$ in $G$. This enables us to state the following theorem. 

\begin{theorem}Let $(H,R)$ be a semi-simple and cosemi-simple triangular
Hopf algebra over $k=\bar{k}$. Then 
$$
BM(k, H, R)=\cases{H^2(G, k^\cdot)\times{\mathbb Z}_2& if
$u\neq1$, $G\cong U\times G/U$;\cr
H^2(G, k^\cdot)& otherwise\cr} 
$$
where $G$ and $U=\langle u \rangle$ are the group data corresponding to $H$. 
\end{theorem} 

\subsection{The multiplication rules}\label{multrules}

So far we have described $BM(k,k[G], R_u)$ with $u\neq1$ as a central
extension of the classical Brauer group $Br(k)$ of the base field
$k$. Our next aim is the explicit description of the multiplication rules
for $BM(k,k[G], R_u)$. 

We shall deal with the cases $G\not\cong U\times G/U$ and
$G\cong U\times G/U$ separately.
 
\begin{itemize}
\item Let $G\not\cong U\times G/U$. In this case each element in
  $BM(k,k[G], R_u)$
can be written as a product of a central simple algebra with trivial
  $G$-action and an algebra of the form $A^\sigma$. We will first
  describe how to do so. Let $[A]\in BM(k,k[G], R_u)$,
let $\zeta$ be as in Theorem \ref{notdirect} and let
  $\zeta([A])=\overline{\sigma}$. Then $[A\sharp
 \overline{A^{\sigma}}]\in\Ker(\zeta)=Br(k)$ and we can
represent  $[A\sharp \overline{A^\sigma}]$ by
  $(A\sharp\overline{A^\sigma})_{triv}$ with the same underlying algebra and
trivial $G$-action. Therefore,
  $[A]=[(A\sharp\overline{A^\sigma})_{triv}][A^\sigma]$ and we need only to describe the classical
  Brauer group class of $(A\sharp\overline{A^\sigma})_{triv}$. 
We recall that $A^\sigma$ is a
  matrix algebra and it is graded central simple, so the same holds
  for its graded opposite. Let us observe that if conjugation by $f_u$
  induces the $u$-action on $A^\sigma$, then the $u$-action on its
  graded opposite can again be induced by $f_u$, and
  $f_u^2=\sigma(u,u)$ in both cases.  
 If $\sigma(u,u)$ is a square, i.e., if the $U$-action is
 strongly inner, then $A\sharp
 \overline{A^\sigma}\cong A\otimes
 \overline{A^\sigma}$ so the Brauer group class of this element is the
  class of $A$. If $\sigma(u,u)$ is not a square 
 \cite[Theorem
  IV.3.8 part (4)]{lam} for $D=F$ implies that $A^\sigma$ is
  isomorphic to the
  tensor product of a matrix algebra with trivial ${\mathbb
  Z}_2$-grading and the quaternion algebra $\langle\frac{-\sigma(u,u),
  1}{k}\rangle$ generated by odd $x$ and $y$ subject to the relations
  $x^2=-\sigma(u,u)$, $y^2=1$ and $xy+yx=0$. The graded
  opposite
  $\overline{A^\sigma}$ is thus isomorphic to the
  tensor product of a matrix algebra with trivial ${\mathbb
  Z}_2$-grading and the quaternion algebra 
$\langle\frac{\sigma(u,u),
 -1}{k}\rangle$, generated by odd $x$ and $y$ subject to the relations
  $x^2=\sigma(u,u)$, $y^2=-1$ and $xy+yx=0$. When we refer to an
  algebra $B$ with trivial $G$-action, we shall denote it by $B_{triv}$. 

Summarizing we have:
$$[A]=\cases{[A_{triv}][A^\sigma]& if $\sigma(u,u)\in(k^\cdot)^2$\cr
[(A\sharp\langle\frac{\sigma(u,u),
 -1}{k}\rangle)_{triv}][A^\sigma]& if $\sigma(u,u)\not\in(k^\cdot)^2.$\cr}$$  
Next we would like to determine the product of these elements. The
only case needing an analysis is the product $[A^\sigma][A^\omega]$.
We have:
$$
\begin{array}{rl}
[A^\sigma\sharp A^\omega]&=[(A^\sigma\sharp A^\omega\sharp\langle\frac{\sigma*\omega(u,u),
 -1}{k}\rangle)_{triv}][A^{\sigma\sharp\omega}]\\
&=[(\langle\frac{-\sigma(u,u),
 1}{k}\rangle\sharp\langle\frac{-\omega(u,u),
 1}{k}\rangle\sharp\langle\frac{\sigma*\omega(u,u),
 -1}{k}\rangle)_{triv}][A^{\sigma\sharp\omega}]\\
&=[(\langle\frac{-\sigma(u,u),
 1}{k}\rangle\sharp\langle\frac{-\omega(u,u),\sigma*\omega(u,u)}{k}\rangle)_{triv}][A^{\sigma\sharp\omega}]\\
&=[\langle\frac{\sigma(u,u),\omega(u,u)
 }{k}\rangle_{triv}][A^{\sigma\sharp\omega}]\\
\end{array}
$$ 
where for the last two equalities we have used the computations in \cite[Lemma
  V.3.2]{lam}.    

\item Let $G\cong U\times G/U$. In this case, if $[A]\in BM(k, k[G], R_u)$ then either
  $A$ is a central simple algebra, and it can be treated as in the
  previous case, or $A\cong B\sharp C(1)$, with $B$ central
  simple and trivial $G/U$-action on $C(1)$. In this case we also have 
$[A\sharp C(-1)]=[B\sharp C(1)\sharp C(-1)]=[B]$. Let $A$ be not central simple and let $\zeta$ be
  as in Theorem \ref{sidirect}. Then $\zeta([A])=(\overline{\sigma},
  -1)$ for some $\sigma$ and an argument similar to the one in the
  previous case applied to $A\sharp C(-1)$ gives
 $[A]=[(A\sharp C(-1)\sharp\langle\frac{-\sigma(u,u),
  1}{k}\rangle)_{triv}][A^\sigma][C(1)]$. The multiplication rules
  follow from those in the previous case once we understand the product
  for $[C(1)][C(1)]=[C(1)\sharp C(1)]$. The cocycle $\sigma$ associated with
  this element can be chosen to be trivial everywhere except from 
$\sigma(u,u)=-1$. Then we have
$$
\left[C(1)\sharp
  C(1)\right]=\left[\left(\left\langle\frac{1,1}{k}\right\rangle\sharp\left\langle\frac{-1,-1}{k}\right\rangle\right)_{triv}\right][A^\sigma]=[A^\sigma].$$
The reader is invited to compare these rules with those in \cite[Theorem V.3.9]{lam}, being
alert that not all representatives in this paper have been chosen as
in \cite{lam}. For instance, the computation of the square of the
triple in \cite{lam}
corresponding to $C(1)$ yields $\langle\frac{1,1}{k}\rangle_{triv}\sharp
\langle\frac{1,1}{k}\rangle$. The first term is a matrix algebra with
trivial action while
the second term lies in the same class as our $A^\sigma$.
\end{itemize}

\section{The Brauer group of $(k[G]\ltimes\wedge V, R_u)$}

Let $(H, R)=(k[G]\ltimes\wedge V,R_u)$ be as in Section \ref{supergroup}.
%
%
Let us denote
by $\rho(g)$ the matrix performing the action of the element $g$ by
conjugation on $V$. Let $v_1,\ldots,v_n$ denote a basis for $V$.
We choose to write $g.v_i=\sum\rho(g)_{ji}v_j$ so
that $\rho$ is a group morphism. 

The triangular Hopf subalgebra $(k[G], R_u)$ is also a quotient of $H$. The
corresponding inclusion and projection maps
$\iota$ and $\pi$ are both
quasitriangular. Therefore, the pull-back $\iota^*$ along $\iota$ induces a
surjective and split map $BM(k, H, R_u)\to BM(k, k[G], R_u)$ which we
still denote by $\iota^*$. Since the Brauer group of a triangular Hopf algebra
is abelian, $BM(k, H, R_u)\cong BM(k, k[G], R_u)\times {\rm
Ker}(\iota^*)$.
We shall devote the rest of this Section to the computation of ${\rm
Ker}(\iota^*)$.

The Hopf subalgebra generated by $u$ and by the
$(u,1)$-primitive
elements in $V$ is
isomorphic, as a triangular Hopf algebra, to the Hopf algebra $E(n)$ 
in Example \ref{enne}
with $R$-matrix $R_u=R_0$. The
corresponding inclusion map $i_E$ is a quasitriangular map and the
pull-back along 
 $i_E$ defines a group morphism $BM(k, H, R_u)\to BM(k, E(n),
 R_u)$ which we shall denote by $i^*_E$. The group $BM(k, E(n),
 R_u)$ has been computed in \cite{GioJuan3} through an analysis of the 
 kernel of the split group epimorphism $j^* \colon BM(k, E(n),
 R_u)\to BW(k)$ induced by the restriction of the $E(n)$-action to a
 $U$-action. One has an isomorphism $\chi\colon \Ker(j^*)\to S^2(V^*)$ so
 $BM(k, E(n), R_u)\cong BW(k)\times S^2(V^*)$.

The restriction of the morphism $i^*_E$ to ${\rm Ker}(\iota^*)$
has image in the Kernel of $j^*$ because $[A]\in{\rm Ker}(\iota^*)$ if
and only if $A$ is an endomorphism algebra with 
strongly inner $G$-action
and $[B]\in{\rm Ker}(j^*)$ if and only if 
$B$ is an endomorphism algebra with 
strongly inner $U$-action.
Thus we have the following commutative
diagram.
$$
\CD
1\to \Ker(\iota^*)@>>>BM(k, H, R_u)@>{\iota^*}>>BM(k, k[G],
R_u)\to 1\\
 @VV{i^*_E}V  @VV{i_E^*}V @VV i^*V\\
1\to\Ker(j^*)@>>>BM(k,E(n), R_u)@>{j^*}>>BW(k)\to1\\
@V\cong V \chi V @. @.\\
S^2(V^*) @.
\endCD
$$
We shall denote the composition $\chi\circ i^*_E\colon {\rm 
Ker}(\iota^*)\to S^2(V^*)$ by $\delta$.

Let us recall how the map $\chi\colon{\rm Ker}(j^*)\to
S^2(V^*)$ was defined in \cite{GioJuan3}. If $A$ represents an element in $\Ker(j^*)$, it is an
endomorphism algebra. The $E(n)$-action is thus inner and the
$U$-action is strongly inner. The subalgebra of $A$ 
generated by the images of $E(n)$ under a map $f$ realizing the action, 
satisfies the relations $f(h)f(l)=\sum f(\h1\l1)c(\h2,\l2)$
for some 2-cocycle $c$. We can make sure that $c$ is lazy. 
The map associating a lazy $c$ to $[A]$ is
well-defined up to (lazy) coboundaries. 
We can always find a privileged representative $\sigma$ 
of the lazy cohomology class of $c$
such that $\sigma(u,u)=1$, $\sigma(u,v_i)=\sigma(v_i,u)=0$ and   
 $\sigma(v_i,v_j)=\Sigma_{ij}$ is a symmetric matrix. In this fashion,
$\overline\sigma$ 
determines an element $\Sigma$ in $S^2(V^*)$. It is shown in
\cite{GioJuan3} that $\chi$ is a group isomorphism.      

\medskip

If $[A]\in {\rm Ker}(\iota^*)$ the $H$-action on $A$ is inner.
If $f$ denotes the realizing map $f\colon H\to A$, 
centrality of $A$ implies that 
$f(h)f(l)=\sum f(\h1\l1)\sigma(\h2,\l2)$ for some    
right 2-cocycle $\sigma$. 
In particular, the restriction of $\sigma$ to
$k[G]\otimes k[G]$ is a coboundary. 

\begin{lemma}\label{sure}Let $[A]\in\Ker(\iota^*)$. We can always find 
a map realizing the action so that its
  restriction to $k[G]$ is an algebra morphism and for
  the  corresponding cocycle $\sigma$ there holds $\sigma(g,
  hv)=\sigma(hv,g)=0$ for every $g, h\in G$ and every $v\in V$.  
\end{lemma} 
\pf Let us choose $f$ so that its restriction to $E(n)$ gives the
privileged cocycle chosen in \cite{GioJuan3}. Let us denote the
corresponding cocycle on $H$ by $\omega$. The restriction of $\omega$ to
$k[G]\otimes k[G]$ is a coboundary $\partial(\gamma)$. Let us extend
$\gamma$ trivially outside $k[G]$. Replacing $f$ by $f*\gamma^{-1}$ 
changes $\omega$ into $\omega^\gamma=(\gamma^{-1}\circ
m)*\omega*(\gamma\otimes\gamma)$ and by construction
$\omega^\gamma(g,h)=1$ for every $g,h\in G$. For $g=h=u$ this shows
that $\gamma(u)=\pm1$. The cocycle $\omega^\gamma$ satisfies:
$$\omega^\gamma(u,u)=1;\quad
\omega^\gamma(u,v_i)=\omega^\gamma(v_i,u)=0;\quad \omega^\gamma(u^av_i,u^bv_j)=\omega(u^av_i,u^bv_j).$$ 
If $f*\gamma^{-1}$ does not satisfy the second condition, we replace $f*\gamma^{-1}$ by $\bar f=f*\gamma^{-1}*\mu$ with $\mu\colon H\to k$
given by $\mu=\varepsilon+\sum_i\omega^\gamma(g,v_i)(gv_i)^*$. 
Then $\bar f(g)\bar f(v)=\bar f(gv)$ and the new cocycle $\sigma$ is such that
$\sigma(g, v)=0$ for every $g\in G$ and $v\in V$. Since
$\mu=\varepsilon$ on $k[G]$ the restriction of $\bar f$ to $k[G]$
coincides with $f*\gamma^{-1}$ which is an algebra morphism.
Besides, we have $\bar f(v)\bar f(u)=\bar f(vu)$, as it was for $f$
and for $f*\gamma^{-1}$ because $\omega(u,v_i)=0$.
%
%
The right cocycle condition gives, for every $g,h\in G$ and every
$v\in V$:
$$\sigma(g,hv)\sigma(h,1)+\sigma(g,hu)\sigma(h,v)=
\sigma(g,h)\sigma(gh,v)$$
hence $\sigma(g,hv)=0$. It follows that for every $h,g\in G$,  
$$
\begin{array}{rl}
\bar f(h)\bar f(v_i)\bar f(g)&=\bar f(hv_i)\bar f(g)\\
&=\bar f(hv_ig)\sigma(h,g)+\bar f(hug)\sigma(hv_i, g)\\
&=\bar f(hg g^{-1}v_ig)+\bar f(hug)\sigma(hv_i, g)\\
&=\sum \rho(g^{-1})_{ji}\bar f(hg v_j)+\bar f(hug)\sigma(hv_i, g)\\
&=\sum \rho(g^{-1})_{ji}\bar f(hg)\bar f(v_j)+\bar f(hug)\sigma(hv_i, g).
\end{array}
$$
We observe that $\bar f(v)$ and $\bar f(u)$ skew-commute because $u$
and $v$ do so. Hence, the first term skew-commutes with $\bar f(u)$ 
and so does each summand in the last term of
the chain of equalities, except from $\bar f(hug)\sigma(hv_i, g)$,
which must be zero. Since $\bar f(hug)$ is invertible, $\sigma(hv_i, g)=0$
and we have the statement.\hfill$\Box$ 

\medskip

Let us observe that the above Lemma is a generalization of the
procedure \cite[Lemma 3.3]{GioJuan3}. In particular, it is
compatible with the choice of the realizing map used for the
construction of $\chi$ because the applied modifications do not change
the value of the cocycle on the pairs $(v_i,v_j)$. 

\medskip

%
%
\begin{lemma}\label{injective}The group morphism $\delta$ is injective. 
\end{lemma}
\pf Since $\delta=\chi\circ i^*_E$  and $\chi$ is injective,
${\rm Ker}(\delta)={\rm Ker}(\iota^*)\cap{\rm Ker}(i^*_E)$. Hence, if
$A$ represents an element in ${\rm Ker}(\delta)$ and $f$ is a
realizing map for the action on $A$, we
may always choose $f$ so that its restriction to $k[G]$  
is an algebra morphism and such that $\sigma(u,v_i)=\sigma(v_i,u)=0$
for every $i$. Since $[A]\in\Ker(i_E^*)$, the restriction of $\sigma$ to 
$E(n)\times E(n)$ is a right coboundary $\mu\circ
m*(\mu^{-1}\otimes\mu^{-1})$ for some cochain $\mu$. Since
$\sigma(u,u)=1$ we have $\mu(u)=\pm1$. 
We can always make sure that $\mu(u)=1$. Indeed, if $\mu(u)=-1$, we can
replace $\mu$ with $\mu*(1^*-u^*)$. This is possible because 
$1^*-u^*$ is a lazy cochain and an
algebra morphism and, by \cite[Lemma 1.6]{lazy},
$(\partial(\mu*(1^*-u^*)))^{-1}=(\partial(\mu))^{-1}$.
Let us extend $\mu$ by $\varepsilon$ outside
$E(n)$. If we replace $f$ by $\bar f=f*\mu$, the new realizing map $\bar f$ 
is an algebra morphism when restricted to both $k[G]$ and $E(n)$.  
To the cocycle $\sigma^\mu$ we can apply the procedure used in Lemma
\ref{sure} to obtain $\sigma^\mu(g,hv)=\sigma^\mu(gv,h)=0$. Such a
procedure can be iterated for increasing 
$l=|P|$ using 
$\gamma_l:=\varepsilon+\sum_{g\in G, |P|=l}\sigma^{l-1}(g, v_P)(gv_P)^*$.
At each step we do not modify the realizing map on the
part of the filtration of $H$ given by $\oplus_{j\le t}k[G]\ltimes\wedge^t V$
corresponding to $t<l$. Therefore the $l$-th step 
ensures that $\sigma(g,hv_P)=0$ for $|P|\le l$. Iterating it up to
$l=n$ ensures that there exists a realizing map $\tilde{f}$ for which
$\sigma(hv_i, g)=0$ and 
$\tilde{f}(gv_P)=\tilde{f}(g)\tilde{f}(v_P)$ for every $g, h\in G$, every $i$ and every $P$.
Besides, the restriction of such an $\tilde{f}$ to $k[G]$ and $E(n)$ is still
a morphism thus $\tilde{f}$ is defined on $H$ exactly by
$\tilde{f}(gv_P)=\tilde{f}(g)\tilde{f}(v_{p_1})\cdots \tilde{f}(v_{p_l})$.
Since the defining 
algebra relations of $H$ are either those in $E(n)$, those in $k[G]$
or those determined 
by $v_ig=g(g^{-1}v_ig)=\sum_{j}\rho(g^{-1})_{ji}gv_j$, all of them are
preserved by $\tilde{f}$ because
$$
\tilde{f}(v_i)\tilde{f}(g)=\tilde{f}(v_ig)=\sum_{j}\rho(g^{-1})_{ji}
\tilde{f}(gv_j)=\sum_{j}\rho(g^{-1})_{ji} \tilde{f}(g) \tilde{f}(v_j)
$$ 
where for the first equality we have used that $\sigma(v_i,g)=0$. 
Therefore $\tilde{f}$ is an algebra morphism, $A$ is a matrix algebra with
strongly inner $H$-action, $[A]=1$ and $\delta$ is
injective.\hfill$\Box$   

\begin{lemma}\label{invariant}The Kernel of $\iota^*$ is isomorphic to
  $S^2(V^*)^G$,  
i.e., to the group of symmetric matrices that are $G$-invariant
with respect to the right action $\Sigma.g=\rho(g)^t\Sigma \rho(g)$ for
every $g\in G$.  
\end{lemma} 
\pf Let $\Sigma$ be the symmetric matrix with coefficients
$\sigma(v_i,v_j)$ determined as in \cite{GioJuan3} and let $f$ be
chosen as in Lemma \ref{sure}, so that
$g.f(v_i)=f(g)f(v_i)f(g)^{-1}=f(gv_ig^{-1})=f(g.v_i)$. Then 
for every $g\in G$ and $i,j\in\{1,\ldots,n\}$ 
$$
\begin{array}{rl}
0&=f(g)(f(v_iv_j+v_jv_i))f(g)^{-1}\\
&=f(g)(f(v_i)f(v_j)+f(v_j)f(v_i)-2\Sigma_{ij})f(g)^{-1}\\
&=(g.f(v_i))(g.f(v_j))+(g.f(v_j))(g.f(v_i))-2\Sigma_{ij}\\
&=(f(g.v_i))(f(g.v_j))+(f(g.v_j))(f(g.v_i))-2\Sigma_{ij}\\
&=f(g.(v_iv_j+v_jv_i))+2\sigma(g.v_i,g.v_j)-2\Sigma_{ij}\\
&=2(\rho(g)^t\Sigma\rho(g))_{ij}-2\Sigma_{ij}
\end{array}
$$
hence $\Sigma\in S^2(V^*)^G$. Conversely, if $\Sigma\in S^2(V^*)^G$,
let us consider $\omega_\Sigma:=r_0*r_{-\Sigma}$ where $r_0$ and
$r_{-\Sigma}$ are as in Example \ref{enne}. By \cite[Example 1.4]{lazy},
$\omega_{\Sigma}$ is a lazy cocycle for $E(n)$. We have, with notation
as in Example \ref{enne}:
$$\omega_\Sigma=\sum_P (-1)^{{|P|(|P|-1)}\over2}\sum_{F, |F|=|P|, \eta\in
S_{|P|}}sign(\eta)(\Sigma_{P,\eta(F)})\bigl((v_P)^*\otimes
(v_F)^*$$
$$+(uv_{P})^*\otimes (v_F)^*+(-1)^{|P|}(v_P)^*\otimes
(uv_F)^*+(-1)^{|P|}(uv_{P})^*\otimes (uv_F)^*\bigr).$$
By direct computation 
\begin{equation}\label{bene}
\omega_\Sigma(u^av_P,u^{b+m}v_Q)=(-1)^{b|P|}\omega_\Sigma(v_P,u^mv_Q)
\end{equation} 
$$\omega_\Sigma(v_P,v_Q)=\cases{0& if $|P|\neq|Q|$\cr
(-1)^{|P|\choose 2}{\rm det}_{PQ}(\Sigma) & if $|P|=|Q|$}$$
where ${\rm det}_{PQ}$ denotes the minor corresponding to the
rows indexed by elements in $P$
and columns indexed by the elements in $Q$.  
A straightforward computation shows that
$$\omega_\Sigma(g.v_P,g.v_Q)=\cases{0& if $|P|\neq|Q|$\cr
(-1)^{|P|\choose 2}{\rm det}_{PQ}(\rho(g)^t\Sigma\rho(g)) & if $|P|=|Q|$}$$
and $G$-invariance of $\Sigma$ implies that
$\omega_\Sigma(g.a,g.b)=\omega_\Sigma(a,b)$ for every $a,b\in E(n)$. 
%
%
Thus, $\omega_\Sigma$ is a $G$-invariant lazy cocycle on $k[U]\ltimes\wedge V\cong
E(n)$ with 
$\omega_\Sigma(v_i,v_j)=\Sigma_{ij}$. 
Let us define: 
\begin{equation}\label{cociclo}\lambda(gv,hw):=\omega_\Sigma(h^{-1}.v,w)
\quad\forall v,w\in\wedge E(n),\ g,h\in G.
\end{equation}

Then $\lambda$ is well-defined thanks to (\ref{bene}), it 
coincides with $\omega_\Sigma$ on $E(n)\otimes E(n)$,
and it is a $2$-cocycle on $k[G]\ltimes\wedge V$. Indeed, for every
$a,b,c\in\wedge V,\ g,h,l\in G$ we have:
$$
\begin{array}{rl}
\sum\lambda((ga)_{(1)},(hb)_{(1)})&\lambda((ga)_{(2)}(hb)_{(2)},lc)\\
&=\sum\omega_\Sigma(h^{-1}.\a1,\b1)\omega_\Sigma(l^{-1}.((h^{-1}.\a2)\b2), c)\\
&=\sum\omega_\Sigma(h^{-1}.\a1,\b1)\omega_\Sigma((h^{-1}.\a2)\b2, l.c)\\
&=\sum\omega_\Sigma(\b1,l.\c1)\omega_\Sigma((h^{-1}.a),\b2 l.\c2)\\
&=\sum\omega_\Sigma(l^{-1}.\b1,\c1)\omega_\Sigma(l^{-1}.(h^{-1}.a),(l^{-1}.\b2)
\c2)\\ 
&=\sum\lambda(\b1,l\c1)\lambda(a,hl(l^{-1}.\b2) \c2)\\
&=\sum\lambda(h\b1,l\c1)\lambda(a,h\b2 l\c2)\\
&=\sum\lambda((hb)_{(1)},(lc)_{(1)})\lambda(ga,(hb)_{(2)}(lc)_{(2)})
\end{array}
$$
hence, since
$\lambda(1,hw)=\omega_\Sigma(1,w)=\varepsilon(w)=\omega_\Sigma(w,1)=\lambda(w,1)$,
$\lambda$ is a 2-cocycle for $k[G]\ltimes\wedge V$. Besides, $\lambda$
is lazy because for every
$a,b\in\wedge V,\ g,h\in G$ we have:
$$
\begin{array}{rl}
\sum\lambda((ga)_{(1)},(hb)_{(1)})(ga)_{(2)}(hb)_{(2)}&=
\sum\omega_\Sigma(h^{-1}.\a1,\b1)gh 
(h^{-1}.\a2)\b2\\  
&=\sum\omega_\Sigma(h^{-1}.\a2,\b2)gh (h^{-1}.\a1)\b1\\
&=\sum\lambda(\a2,h\b2)gh (h^{-1}.\a1)\b1\\
&=\sum\lambda((ga)_{(2)},(hb)_{(2)})(ga)_{(1)}(hb)_{(1)}.
\end{array}
$$
If we construct a $(k[G]\ltimes\wedge V,R_u)$-Azumaya algebra $A^\lambda$
associated to $\lambda$ as in
\cite[Lemma 4.8]{GioJuan3}, then the $G$-action is strongly inner because
$\lambda(g,h)=1$ for every $g,h\in G$, so $[A^\lambda]\in{\rm
Ker}(\iota^*)$. The restriction of $\lambda$ to $(k[U]\ltimes\wedge
V)^{\otimes 2}$ is
a 2-cocycle with $\lambda(u,u)=1$,  $\lambda(u,v_i)=\lambda(v_i,u)=0$ and
$\lambda(v_i,v_j)=\omega_\Sigma(v_i,v_j)=\Sigma_{ij}$ (that is, its
restriction to $E(n)$ is lazy
cohomologous to the privileged $\sigma$ with
$\sigma(v_i,v_j)=\Sigma_{ij}$ in \cite{GioJuan3}). Hence $\delta({\rm
  Ker}(\iota^*))=S^2(V^*)^G$ 
and we have the statement.\hfill$\Box$

%
\begin{remark}{\rm The reader is invited to compare formula
(\ref{cociclo}) with  \cite[Formula
(4.20)]{Juan-Florin}) where a different language for an analogous
construction is used.} 
\end{remark}    

Combining all results in this Section we have the following theorems.

\begin{theorem} Let $(k[G]\ltimes\wedge V, R_u)$ be a modified
  supergroup algebra. Then 
$$BM(k,k[G]\ltimes\wedge V, R_u)\cong BM(k,k[G],R_u)\times S^2(V^*)^G.$$
\end{theorem}

\begin{theorem}\label{tutte}
Let $(H, R)$ be a finite-dimensional triangular Hopf
algebra over $k={\overline k}$ with ${\rm char}(k)=0$. Then  
$$BM(k, H, R)=\cases{H^2(G,
    k^\cdot)\times{\mathbb Z}_2\times S^2(V^*)^G& if $u\neq1$ and $G\cong
    U\times G/U$\cr
H^2(G, k^\cdot)\times S^2(V^*)^G & otherwise\cr}
$$ where $G,\,u,\,U=\langle u\rangle,\,V$ are the corresponding data for $(H, R)$.
\end{theorem}

The multiplication rules in Section \ref{multrules} together with the
relation $[A^\sigma][A^\omega]=[A^{\sigma*\omega}]$ for
representatives $\sigma$ and $\omega$ of elements in $H^2_L(E(n))$, provide the
multiplication rules for $BM(k, k[G]\ltimes\wedge V, R_u)$. In this
case, if $\delta([A^\sigma])=\Sigma$ and $\delta([A^\omega])=\Omega$
then $\delta([A^{\sigma*\omega}])=\Sigma+\Omega$. 

\section{Brauer groups and lazy cohomology}

In this Section we shall apply some arguments  similar to those 
in \cite{Juan-Florin} in
order to compute $H^2_L(H)$ for $H=k[G]\ltimes \wedge V$. We shall avoid the
language of Yetter-Drinfeld modules. Indeed, in this particular case
not so much ado is needed because $\wedge V$ and $E(n)$ are not too
different from each other. When $V=0$ we have $H=k[G]$ and lazy cohomology coincides
with usual group cohomology. We shall assume that $V\neq0$.

Let $\sigma$ be a $2$-cocycle for $H$. Then 
its restrictions $\sigma_E$ and $\sigma_G$ to $E(n)^{\otimes 2}$ and
to $k[G]^{\otimes2}$, respectively, are
both $2$-cocycles. If $\sigma$ is lazy then $\sigma_E$ is lazy. If
$\sigma$ is a coboundary then $\sigma_E$ will be a coboundary, too. 
If $\sigma$ is a lazy $2$-cocycle for $H$, then the lazy condition applied to
$g$ and $hv$ (respectively, $hv$ and $g$) shows that $\sigma_G(g,
hu)=\sigma_G(g,h)$ and $\sigma_G(hu,g)=\sigma_G(h,g)$ for every $g,
h\in G$. Therefore $\sigma_G$ is doubly $U$-invariant and it
determines a cocycle $\sigma_{G/U}$ for $G/U$. Besides, if $\sigma$ is a lazy
coboundary, then $\sigma=\partial(\gamma)$ for some $\gamma$ with
${\rm ad}(\gamma)=\id$. Then ${\rm ad}(\gamma)(gv)=gv$ implies
that 
$\gamma(gu)=\gamma(g)$, so $\sigma_{G/U}$ is a coboundary for
$G/U$. Thus, the assignment $\sigma\mapsto(\sigma_E,\sigma_{G/U})$
induces a group morphism ${\rm Bires}\colon
H^2_L(H)\to H^2_L(E(n))\times H^2(G/U,{\mathbb C}^\cdot)$. The rest of this section will
be devoted to the analysis of this morphism. 

\begin{lemma}\label{kernel}
The morphism ${\rm Bires}$ is injective.
\end{lemma}     
\pf For every lazy $2$-cocycle $\sigma$ for $H$ we consider 
the bicleft extension $k\sharp_\sigma H$
with cleaving map $f$. We shall show that if $\sigma$
represents an element in $\Ker({\rm Bires})$, then we may replace $f$
with an algebra map. This will imply that $\sigma$ is a
coboundary. Then we shall show that we can choose it to be a lazy coboundary. 
Since $\sigma_{G/U}$ is a coboundary $\partial(\gamma)$, 
the cochain which is trivial outside $k[G]$ and defined by
$\gamma\pi$ on $k[G]$, with $\pi$ the natural projection on $G/U$, is
lazy. Hence we might as well replace $\sigma$ by
$\sigma^{(\gamma\pi)^{-1}}=
((\gamma\pi)^{-1}\otimes(\gamma\pi)^{-1})*\sigma*(\gamma\pi)\circ m$ 
without changing its restriction $\sigma_E$ to $E(n)^{\otimes 2}$. 
Therefore we may always assume that the restriction of $\sigma$ to $G$ 
is trivial so that the restriction of the cleaving map $f$ to $k[G]$
is an algebra map. Besides $\sigma_E=\partial(\mu)$ for some lazy cochain
$\mu$. Let us observe that, in particular, $\mu(u)=1$. 
If we extend $\mu$ to $H$ by $\mu(gv_P)=\varepsilon(gv_P)$ if $g\in G$,
$g\not\in U$, then $\mu$ is lazy for $H$. If we replace $\sigma$ by
$\sigma^{\mu^{-1}}=(\mu^{-1}\otimes\mu^{-1})*\sigma*\mu\circ m$ we
have $\sigma_G=\varepsilon\otimes\varepsilon$ 
and $\sigma_E=\varepsilon\otimes\varepsilon$. Hence, we may always
assume that the restrictions of the cleaving map $f$ to $k[G]$ and to
$E(n)$ are algebra maps. Arguments similar to those used to prove 
Lemma \ref{injective} show that 
$\sigma$ is cohomologous to a $2$-cocycle for which $f$ is an algebra
map, i.e., $\sigma\in Z^2_L(H)\cap B^2(H)$. By 
Lemma \ref{cappa}, $\sigma=\partial(\gamma)$ with $\gamma$ either lazy 
or $\ad(\gamma)$ equal to conjugation by $u$. In the latter case, the
proof of Lemma \ref{cappa} shows that $\ad(\gamma)=\ad(\zeta)$ for
some $\zeta$ trivial outside $k[G]$. Therefore $\gamma\circ
\zeta^{-1}$ is lazy  and $\sigma$  
 is lazy cohomologous to $\partial(\zeta)$. By hypothesis ${\rm Bires}(\bar
\sigma)={\rm Bires}(\overline{\partial(\zeta)})=1$ so the
restriction of $\partial(\zeta)$ to $G$ is a coboundary for $G/U$.
Thus, we may choose $\zeta$ to be $U$-invariant on $G$. This
implies that $\zeta$ is lazy for $H$, whence the statement.\hfill$\Box$ 

\smallskip

\begin{remark}{\rm We could also replace $\rm Bires$ by  a morphism
$H^2_L(H)\to S^2(V^*)^G\times {\rm Infl}(H^2(G/U, k^\cdot))$. The
 Kernel of this morphism would be isomorphic to $K(H)$. This should
 explain why $K(H)\cong{\rm Ker}({\rm Infl})$ (see Section \ref{sharp}).}
\end{remark}

\smallskip

\begin{lemma}\label{image}The image of ${\rm Bires}$ is 
$S^2(V^*)^G\times H^2(G/U,k^\cdot)$. 
\end{lemma}  
\pf The subgroup $1\times H^2(G/U,k^\cdot)$ of $S^2(V^*)\times
H^2(G/U,k^\cdot)$ is contained in ${\rm Im}({\rm Bires})$. Indeed, any
cocycle
$\sigma$ on $G/U$ can be extended to a lazy cocycle on $H$
by: $\sigma(gv_P, h v_Q):={\rm Infl}(\sigma)(g,
h)\delta_{P,\emptyset}\delta_{Q,\emptyset}$. 

Besides, it follows from the construction of the cocycle $\lambda$ in 
the proof of Lemma \ref{invariant} that the subgroup 
$S^2(V^*)^G\times 1$ of $S^2(V^*)\times
H^2(G/U,k^\cdot)$ is also contained in ${\rm Im}({\rm Bires})$.
On the other hand, if $\sigma$ is a lazy cocycle for $H$ then
$\sigma_E$ is lazy cohomologous to a $2$-cocycle $c$ with 
$c(v_i,v_j)$ symmetric. If $\gamma$ is the $E(n)$-lazy cochain
for which $\sigma_E^\gamma=c$, then $\gamma(u)=1$ and 
we may extend $\gamma$ on $H$ to a lazy cochain by $\gamma(g)=1$ for
every $g\in G$ and $\gamma(t)=0$ if $t\not\in k[G]\cup
E(n)$. Replacing $\sigma$ by $\sigma^\gamma$ 
we obtain
$\sigma(v_i,v_j)=\Sigma_{ij}$ with $\Sigma$ symmetric.  
The lazyness condition on $g$ and $hv$ ($hv$ and $g$,
respectively) shows that we necessarily have $\sigma(g,
hv_i)=0=\sigma(hv_i, g)$ for every $g, h\in G$ and for every $i$ and   
a computation similar to the one in
the proof of Lemma \ref{invariant}, applied to $k\sharp_\sigma H$  
shows that $\Sigma\in S^2(V^*)^G$.   
Thus 
\begin{equation}
S^2(V^*)^G\times 1\subset {\rm Im}({\rm Bires})\subset
S^2(V^*)^G\times H^2(G/U,k^\cdot). 
\end{equation}
The above inclusions yield the statement. \hfill$\Box$

\smallskip

We have proved the following Theorem.

\begin{theorem}\label{pigra}Let $H=k[G]\ltimes\wedge V$. Then $H^2_L(H)\cong
S^2(V^*)^G\times H^2(G/U,k^\cdot)$. 
\end{theorem}

\smallskip

\begin{remark}{\rm Combining the approach in \cite{Juan-Florin} with
  the above description it should be possible to compute the lazy
  cohomology group for Radford biproducts, at least when the acting
  Hopf algebra is cocommutative. This is the subject of a joint
  forthcoming paper.}
\end{remark}

The computations above show that, even if there is a relation between
lazy cohomology and the Brauer group $BM$, in general lazy cohomology
is not a direct summand of $BM$, nor always a quotient of it. In
fact, even when $G={\mathbb Z}_2$ and $V=0$, $H^2_L(k[G])$ is not a
subgroup of $BM(k,k[G],R_u)\cong BW(k)$, nor of $BQ(k,k[G])\cong
BD(k, k[{\mathbb Z}_2])$, computed in \cite{DMF}.  On the other hand,
the linear summand of lazy cohomology is always a direct summand of
the Brauer group.

\section{Comparison with well-known exact sequences}\label{comparison}

In this Section we compare our results with some known exact sequences
involving Brauer groups and groups of Galois objects. This
comparison involves lenghty but not enlightening computations, so we will
skip the details whenever possible.

\smallskip

Let for the moment $G$ be a finite abelian group. If $k$ has enough
roots of unity then $G\cong G^*$ and
$$BM(k,k[G],R)\cong BC(k,
    k[G]^*, R)\cong BC(k,
    k[G^*], R)\cong B_R(k,G^*),$$ where $B_R(k,G^*)$ denotes the
    Brauer 
group of $G^*$-graded Azumaya
    algebras as in \cite{Chi}, \cite{CGO} and \cite{knus}. The first
    isomorphism is induced by the usual duality functor and the last one
    depends on the fact that a $G^*$-comodule algebra $A$ a
    $G^*$-graded by $x=\sum_{\chi\in G^*} x_\chi\in
    \sum_{\chi\in G^*}A_\chi=A$ if $\rho(x)=\sum_{\chi\in G^*}
    x_\chi\otimes\chi$ and viceversa. 

In \cite{Chi}, \cite{CGO} an  exact sequence
$$
\CD
1\longrightarrow Br(k)\longrightarrow B_R(k, G)@>\pi>> {\rm
  Galz}_R(k, G)
\endCD
$$
where ${\rm  Galz}_R(k, G)$ is a suitable group of $G$-graded Galois
objects is described. The rightmost arrow is surjective when $k$
is nice enough (e.g., when $k$ is a field, which is our case). In \cite{Orz} a cohomological
interpretation of the group of equivalence classes of Galois
extensions with normal basis is given. 
In our terms, the map $B_R(k, G^*)\longrightarrow {\rm
  Galz}_R(k, G^*)$ associates to the class $[A]$ the equivalence class
of the centralizer, in a suitable representative $A$, of the
$G$-invariant subalgebra of $A$. The representative has to be chosen {\em
$G^*$-fully graded}, in the terminology of \cite{CGO}. We are
concerned with the images of our classes $[A^\sigma]$ and $[C(1)]$. It is not hard to verify that $A^\sigma$ is fully
graded and that $\pi([A^\sigma])$ is the class of the twisted group
algebra $k_\sigma[G]$ with associated cocycle $\sigma$, i.e., the
subalgebra generated by $f_h=f(h)$ with $h\in G$. Since $\rho(f_h)=\sum_{g\in
  G}\frac{1}{\sigma(g,g^{-1})}f_g\circ f_h\circ f_{g^{-1}}\otimes
g^*$, a direct computation shows that $f_h$ has degree $\chi_h$ where
$\chi_h(g)=\sigma(g,h)\sigma^{-1}(h,g)$. By \cite[Page 311]{CGO} the right $G^*$-action on
this representative of $\pi([A^\sigma])$ is given by
$(f_h\leftharpoonup
\chi)=\chi(h)R_u(\chi,\chi_h)f_h=\chi(u)^{|h|_\sigma}f_h$. 

\smallbreak

When $R=R_u$ and $G\cong U\times G/U$ the class $[C(1)]$ is an element in $B_R(k,
 G^*)$. Its representative $C(1)$
is not fully graded because $(G/U)^*$ acts trivially on it. Once we
replace $C(1)$ by $C(1)\sharp \End(k[G^*])$ we see that $\pi([C(1)])$ is
the isomorphism class of $k[G^*]\cong k[U^*]\otimes k[(G/U)^*]$ with the usual
${\mathbb Z}_2$-action on $k[U^*]\cong C(1)$ and the regular
$(G/U)^*$-action on $k[(G/U)^*]$. 
The grading of $k[U^*]$ is the non-trivial $\Z2$-grading while
$k[(G/U)^*]$ is trivially graded. 

\smallskip

The exact sequence in \cite{CGO} has been generalized in
\cite{beattie} to the case of commutative cocommutative Hopf algebras with
trivial $R$-matrices and in \cite{ulbrich} to commutative
cocommutative Hopf algebras with $R$ a bipairing on $H^*\otimes
H^*$. The most general case is dealt with in \cite{sequence} where
 an exact sequence
\begin{equation}\label{seqgal}
\CD
1\longrightarrow Br(k)\longrightarrow BC(k, H , R)@>\pi>> {\tt
  Gal}({\cal H}_R)
\endCD
\end{equation}
is constructed for  $(H, R)$ a dual quasitriangular Hopf algebra. Here
${\cal H}_R$ is the braided Hopf algebra of \cite[Theorem
7.4.1]{mabook} and ${\tt Gal}({\cal H}_R)$ is a group of quantum commutative
biGalois objects for $({\cal H}_R)^*$. In \cite{chenzhang} it is shown
that ${\tt
  Gal}({\cal H}_R)$ is invariant under cocycle twist. The sequence
(\ref{seqgal}) has relevant
theoretical meaning but ${\tt Gal}({\cal H}_R)$ might be very hard to compute. It
is possible that the analysis of the relations between (\ref{seqgal})
and our computations in the triangular case will give an indication
on how to handle $BM(k,H, R)$ in the general quasitriangular case. 

If $(H, R)$ is dual triangular and, for instance, ${\overline k}=k$, then 
$(H, R)$ is the Doi twist of $((k[G]\ltimes\wedge V)^*, R_u)$  for
some $G, V, u$. Zhang's sequence becomes: 
\begin{equation}
\CD
1\longrightarrow Br(k)\longrightarrow BM(k, k[G]\ltimes\wedge V, R_u)@>\pi\psi>> {\tt
  Gal}({\cal H}_R).
\endCD
\end{equation}
Here we have used that $[A]\mapsto[A^{\rm op}]$ defines a group
isomorphism 
$$\psi\colon BM(k,k[G]\ltimes\wedge V,R_u)\to BC(k,(k[G]\ltimes\wedge
V)^*, R_u)$$ stemming from the monoidal functor 
$$({\cal D}, \tau_{UV},
\id)\colon (_{H^*}\!{\cal M},\otimes,k,\id,\id,\id)\to({\cal M}^{H},
\otimes^{\rm rev}, k,\id,\id,\id)$$ where $\cal D$ is the usual duality
functor from the category ${_{H^*}\!\cal M}$ of left $H^*$-modules to
the category ${\cal M}^{H}$ of right $H$-comodules and $\otimes^{\rm rev}$
denotes the opposite tensor product. If $(H^*, R)=(k[G]\ltimes\wedge V, R_u)$, by \cite[Theorem
7.4.2]{mabook} we have ${\cal H}_R^*\cong k[G]\ltimes\wedge V$ as an
algebra, with
the (super cocommutative) Hopf superalgebra structure given on generators 
by $\Delta(g)=g\otimes g$ and
$\Delta(v)=v\otimes 1+1\otimes v$. We shall denote this Hopf
superalgebra by $\cal A$. An element in ${\tt Gal}({\cal H}_R)$ is determined by an
$H$-Yetter-Drinfeld module structure and an ${\cal A}$-bicomodule
structure. The unit element in ${\tt Gal}({\cal
  H}_R)$ is represented by $H^*$, with actions and coactions defined
by \cite[Formulas (3), (6), (11)]{sequence}.  

The map $\pi$ is given as follows. For a
class in $BC(k, H, R_u)$ we consider a representative $A$ which is
Galois. This is always possible by \cite[Corollary
  4.2]{sequence}. Then, as an algebra, $\pi([A])=C_A(A_0)$, the
centralizer of the $H$-coinvariants. The Yetter-Drinfeld module
structure on $C_A(A_0)$
is given by the restrictions of the Miyashita-Ulbrich-Van Oystaeyen (MUVO) action $\rightharpoonup$ (see
\cite[\S 2.3]{CVZ}) of $H$ and of the $H$-coaction $\rho$ on $A$.
The $\cal A$-bicomodule
structure on $C_A(A_0)$ is given by $\rho_1$ and $\rho_2$ which are dual to actions
$-\!\triangleright\;$ and $\triangleleft\!-$. These actions are  
modifications through the
$R$-matrix (\cite[Formulas (3),(6)]{sequence}) and through $\rho$ of
the $H$-action $\rightharpoonup$.

\smallskip

We aim at the description of the image of the classes $[A^\sigma]$
when $H^*=k[G]$ and $R=1\otimes 1, R_u$ and when
$H^*=k[G]\ltimes\wedge V$ and $R=R_u$, and the image of the class
$[C(1)]$ when $H^*=k[G]$ or $H^*=k[G]\ltimes\wedge V$; $R=R_u$ and $G\cong U\times G/U$.

Let us first assume that $V=0$. Then ${\cal A}=H^*=k[G]$ is a genuine
cocommutative  Hopf algebra with $G$ not necessarily abelian. Its dual $k[G]^*$
is the span of $v_h=h^*$ for $h\in G$ with product $v_hv_g=\delta_{h,g}v_h$.  
The unit element in ${\tt
  Gal}({\cal H}_R)$ is represented by $k[G]$ with the regular
left and right coactions $\rho_1$ and $\rho_2$ for $\cal A$;
$k[G]^*$-action given by $v_h.g=\delta_{h,g}g$ and $k[G]^*$-coaction
given by $\rho(g)=\sum_{h\in G}hgh^{-1}\otimes v_h$.    
%

\smallskip

If $R=1\otimes 1, R_u$ then the set $H^2(G, k^\cdot)$ occurs in $BM(k,
k[G], R)$. The image of $[A^\sigma]$ can be computed
as follows. The class of $A^\sigma$ (with product $\circ$) in $BM(k,k[G],R)$ is
mapped by $\psi$ to the class of $[(A^\sigma)^{\rm op}]$. The coinvariants in
$(A^{\sigma})^{\rm op}$ 
correspond to the centralizer of the induced subalgebra, that is, of
the algebra generated by the $f_g$'s in $[A^\sigma]$. As an algebra, $\pi\psi([A^\sigma])$ 
is the equivalence class of the algebra $(k_\sigma [G])^{\rm op}$ with
product $\bullet$.
The $k[G]^*$-comodule structure is given, for both choices of $R$, by: 
$$\rho(f_h)=\sum_{g\in G}\frac{\sigma(g,h)\sigma(gh,
  g^{-1})}{\sigma(g,g^{-1})}f_{ghg^{-1}}\otimes v_g.$$
The MUVO action is, for both choices of $R$, as follows. 
If $\beta$ denotes the natural map $A\otimes_{A_0} A\to A\otimes k[G]^*$ 
and if $(1\otimes v_h)=\beta(\sum F_i(h)\otimes_{A_0} f_i(h))$, then 
$$
\begin{array}{rl}
\delta_{gh}&=\sum F_i(h)\bullet f_{i0}(h)\langle f_{i1}(h),
g\rangle\\
&=\sum F_i(h)\bullet (g.f_i(h))\\
&=\sum\frac{1}{\sigma(g, g^{-1})} F_i(h)\bullet (f_g\circ f_i(h)\circ
f_{g^{-1}})\\
&=\sum\frac{1}{\sigma(g, g^{-1})} F_i(h)\bullet f_{g^{-1}}\bullet f_i(h)\bullet
f_{g}\\
&=(v_h\rightharpoonup f_{g^{-1}})\bullet f_{g^{-1}}^{-1}.
\end{array}
$$
Thus $v_h\rightharpoonup f_{g}=\delta_{h, g^{-1}}f_g$. 

If $R=1\otimes 1$ then the $\cal A$-actions $-\!\triangleright\;$ and
$\triangleleft\!-$ coincide with the $\rightharpoonup$ action. 
The $\cal A$-comodule structures are:
$$\rho_1(f_g)=f_g\otimes
g^{-1};\quad\quad\rho_2(f_g)=g^{-1}\otimes f_g.
$$ 
If $R=R_u$ with $u\neq1$ 
then 
$$v_h-\!\triangleright\;f_g=f_g\triangleleft\!-
v_h=\cases{\delta_{h,g^{-1}}f_g& if $|g|_\sigma=0$;\cr
\delta_{h,ug^{-1}}f_g& if $|g|_\sigma=1$.\cr}$$

The $\cal A$-comodule structures are:
$$\rho_1(f_g)=\cases{f_g\otimes
g^{-1}& if $|g|_\sigma=0$;\cr
f_g\otimes
ug^{-1}& if $|g|_\sigma=1$;\cr} \quad
\rho_2(f_g)=\cases{g^{-1}\otimes f_g & if $|g|_\sigma=0$;\cr
ug^{-1}\otimes f_g & if $|g|_\sigma=1$.\cr}
$$
The assignment $f_g\mapsto f'_{g^{-1}}$ determines an algebra
isomorphism $k_\sigma[G]^{\rm op}\cong k_{\sigma'}[G]$ with
$\sigma'(g,h)=\sigma(h^{-1}, g^{-1})$. Therefore we may write:
$\pi([A^\sigma])\cong k_{\sigma'}[G]$ as an algebra, with: comodule structure
$$\rho(f'_h)=\sum_{g\in G}\frac{\sigma'(g, g^{-1})}{\sigma'(h, g)\sigma'(g,
  hg^{-1})}f'_{ghg^{-1}}\otimes v_g;$$
action $v_h\rightharpoonup f'_g=\delta_{hg}f'_g$; and $\cal A$-comodule
  structures
$$\rho_1(f'_g)=f'_g\otimes
g;\quad\quad\rho_2(f'_g)=g\otimes f'_g
$$ if $R=1\otimes 1$; and 
$$\rho_1(f'_g)=\cases{f'_g\otimes
g& if $|g|_\sigma=0$;\cr
f'_g\otimes
ug& if $|g|_\sigma=1$;\cr} \quad
\rho_2(f'_g)=\cases{g\otimes f'_g & if $|g|_\sigma=0$;\cr
ug\otimes f'_g & if $|g|_\sigma=1$\cr}
$$
if $R=R_u$ with $u\neq1$. 
%
In both cases, if $k^\cdot=(k^\cdot)^2$ then by \cite[Lemma 3.6]{kar2}
we may choose $\sigma$ in its cohomology class so that $\sigma'=\sigma^{-1}$. 
%

\smallskip

When $R=R_u$ and $G\cong U\times G/U$, we analyze the image
of $[C(1)]$. If we take
$A=C(1)\sharp\End(k[G]^*)$ as a representative of this class in $BC(k,
k[G]^*, R_u)$, its image
$\pi([A])$ is as follows. Let $g\mapsto f_g$ for $g\in G$
denote the map realizing the (strongly inner) $k[G]$-algebra action on 
$\End(k[G]^*)^{\rm op}={\cal D}^{-1}(\End(k[G]^*))$. 

The algebra $\pi([A])=C_A(A_0)$ is generated by the elements $1\sharp f_h$ for
$h\in G/U$ and $x\sharp f_u$. Since the usual isomorphism
$C(1)\sharp\End(k[G]^*)^{op}\cong C(1)\otimes\End(k[G]^*)^{op}$ maps
$x\sharp f_u$ to $x\otimes 1$,
the element $x\sharp f_u$ is central in $A$. Then, 
by \cite[Lemma 2.3.1 (a)]{CVZ}  
the MUVO action on it is trivial. Let now $v_g=g^*\in k[G]^*$.
Then the $k[G]^*$-comodule structure is given by:
$$\rho(x\sharp f_u)=(x\sharp f_u)\otimes\sum_{h\in G/U}(v_h-v_{hu})$$
$$\rho(1\sharp f_h)=\sum_{t\in G}(1\sharp f_{tgt^{-1}})\otimes v_t.$$
A computation similar to the case of $[A^\sigma]$ shows that 
$v_h\rightharpoonup (1\sharp f_g)=\delta_{h,g^{-1}}(1\sharp f_g)$ for
$g\in G/U$.
%
%
An
analysis of the structure of the actions $-\!\triangleright\;$ and
$\,\triangleleft\!-$ in \cite[Formulas (3), (6)]{sequence}  
on the elements of $k[G/U]$  shows that they both coincide with the
MUVO-action, while 
$$(x\sharp
f_u)\,\triangleleft\!-v_h=\delta_{h,u}(x\sharp
f_u)=v_h-\!\triangleright (x\sharp f_u).$$
Thus, the corresponding comodule maps are
$$
\begin{array}{ll}
\rho_1(x\sharp f_u)=(x\sharp f_u)\otimes u;&\rho_2(x\sharp
f_u)=u\otimes (x\sharp f_u);\\
\rho_1(1\sharp f_h)=(1\sharp f_h)\otimes h^{-1};&\rho_2(1\sharp
f_h)=h^{-1}\otimes (1\sharp f_h).\\
\end{array}
$$
The assignment $x\sharp f_u\mapsto u$ and $1\sharp f_{h}\mapsto
h^{-1}$ for $h\in G/U$ determines an
algebra isomorphism $\pi\psi([A])\to k[G]$.
Through this isomorphism $\pi\psi([A])\cong k[G]$ with coregular $k[G]$-comodule
structures $\rho_1$ and $\rho_2$ and with $k[G]^*$-comodule structure
$\rho$ given by:
$$\rho(u)=u\otimes\sum_{h\in G/U}(v_h-v_{hu}),\quad
\rho(h)=\sum_{g\in G}ghg^{-1}\otimes v_g$$
for $h\in G/U$. The $k[G]^*$-action is given by $v_h\rightharpoonup g=\delta_{h,g}g$
for $g\in G/U$ and 
$v_h\rightharpoonup u=\delta_{1,h}u$, for every $h\in G$. 

%

\smallskip

Let us now assume that $V\neq0$. Then necessarily $R=R_u$ with
$u\neq1$ and ${\cal A}=k[G]\ltimes\wedge V$ is a Hopf
superalgebra. 
 
We are concerned with the image of $[A^\sigma]$. We shall separately 
deal with the cases: $\sigma$ is lazy, trivial on $k[G]$; $\sigma$ is any cocycle for $G$, trivial
outside $k[G]$. In the first case $A^\sigma$ is isomorphic to
$\End(H^*)$ as an algebra and $\pi\psi([A^\sigma])$ is the
opposite of the induced subalgebra of $A^\sigma$, i.e., the algebra
generated by the $f_h$'s as before, for $h\in H^*$. Let us denote by
$\circ$ the product in $A^\sigma$ and by $\bullet$ its opposite
product. The $H$-comodule structure on $\pi\psi(A^\sigma)$ is given by 
$$\rho(f_h)=\sum_{a\in H^*}f_{a_{(1)}}\circ f_h\circ f_{a_{(2)}}^{-1}\otimes
a^*$$
where the expression $f_{a}^{-1}$ stands for the convolution
inverse of $f_{a}$. It can be proved as in the previous cases that 
the MUVO action on $f_g$ for $g\in G$ is given by $h^*\rightharpoonup
f_g=\delta_{h, g^{-1}} f_g$, for $h^*\in H$ dual to $h$. For the elements $f_v$ with $v\in V$ we
have, if $\beta(\sum F_i(a)\otimes f_i(a))=1\otimes a^*$ for $a^*\in
H$:
$$
\begin{array}{rl}
\delta_{a, uv}&=\sum F_i(a)\bullet f_{i0}(a)\langle f_{i1}(a), uv\rangle\\
&=\sum F_i(a)\bullet (uv.f_i(a))\\
&=\sum F_i(a)\bullet (f_{uv}\circ f_{i}(a)\circ f_u^{-1})+\sum
F_i(a)\bullet (f_i(a)\circ f_{uv}^{-1})\\
&=\sum F_i(a)\bullet f_u\bullet f_i(a)\bullet f_{uv}+
\sum F_i(a)\bullet f_v\bullet f_i(a)\\
&=(a^*\rightharpoonup f_u)\bullet
f_{uv}+a^*\rightharpoonup f_v\\
&=-\delta_{a,u}f_v+a^*\rightharpoonup f_v
\end{array}
$$
where for the fourth equality we have used the formula for $f^{-1}_h$
in \cite[Lemma 4.8]{GioJuan3} and for the last equality we have used
our knowledge of the MUVO action on $f_g$ for $g\in G$.  
Thus, $a^*\rightharpoonup f_v=\delta_{a,
  uv}+\delta_{a,u}f_v$. A direct computation using that: $\Delta(a^*)$
is a linear combination of elements $r^*\otimes s^*$ with $rs=a$ and that
$S$ preserves the filtration induced by powers of $V$  
shows that for every $a\in H^*$, every $g\in G$, and every $v\in V$:
$$
\begin{array}{lll}
a^*-\!\triangleright\;f_g&=f_g\,\triangleleft\!-a^*&=\delta_{g^{-1},a}f_g;\\
a^*-\!\triangleright\;f_v&=f_v\,\triangleleft\!-a^*&=\delta_{a,uv}+\delta_{a,1}f_v.
\end{array}$$
The $\cal A$-bicomodule structure is, for $g\in G$ and $v\in V$:
$$
\begin{array}{ll}
\rho_1(f_g)=f_g\otimes g^{-1};&\rho_1(f_v)=1\otimes uv+f_v\otimes 1;\\
\rho_2(f_g)=g^{-1}\otimes f_g;&\rho_2(f_v)=uv\otimes 1+1\otimes f_v.\\
\end{array}
$$
We observe that $\pi\psi([A^\sigma])$ is isomorphic, as an algebra, to
$(k\sharp_\sigma H^*)^{\rm op}$. The assignment $f_h\mapsto f'_{S^{-1}(h)}$
determines an algebra isomorphism $(k\sharp_\sigma H^*)^{\rm op}\to
H^*\,_{\sigma'}\sharp k$, with right cocycle
$\sigma'$ given by
$\sigma'(a,b)=\sigma(S(a), S(b))$. Since $\sigma$ is lazy,
$\sigma'$ is lazy as well. Through this identification,
$\pi\psi([A^\sigma])$ is $k\sharp_{\sigma'} H^*$ with Yetter-Drinfeld
module structure
$$
\begin{array}{l}
\rho(f'_h)=\sum_{a\in
  H^*}\sigma^{-1}(S\a3,\a4)f_{a_{(2)}}f'_hf'_{S^{-1}(a_{(1)})}\otimes a^*;\\
a^*\rightharpoonup f'_g=\delta_{a,g}f'_g;\quad\quad a^*\rightharpoonup
  f'_{uv}=\delta_{a,uv}+\delta_{a,1}f'_{uv}
\end{array}
$$
and with $\cal A$-bicomodule structure: 
$$
\begin{array}{ll}
\rho_1(f'_g)=f_g\otimes g;&\rho_1(f'_{uv})=f'_{uv}\otimes1+1\otimes
uv;\\
\rho_2(f'_g)=g\otimes f'_g;&\rho_2(f'_{uv})=uv\otimes 1+1\otimes
f'_{uv}.
\end{array}
$$
Let us observe that if $\sigma$ is the privileged cocycle in
\cite{GioJuan3} corresponding to the symmetric matrix $\Sigma$ then
$\sigma'$ is the
privileged cocycle corresponding to $-\Sigma$, i.e., it represents the
inverse cohomology class.

Let now $\sigma$ be a cocycle for $k[G]$ and trivial elsewhere. Since $(A^\sigma)^{\rm op}=(\End(k[G]))^{\rm op}$ is not
Galois, we
consider the representative $(A^\sigma)^{\rm op}\sharp\End(H^{\rm op})$ where
$\End(H^{\rm op})$ is as in \cite[Corollary 4.2]{sequence}. This
corresponds to the representative 
$A\cong ((A^\sigma)^{\rm op}\sharp\End(H^{\rm op}))^{\rm op}$
of the class $[A^\sigma]$ in $BM$.
The usual flip map $\tau$ determines an isomorphism 
$A\cong \End(H^{\rm op})^{\rm op}\sharp A^\sigma$, and $\End(H^{\rm op})^{\rm op}$ has a strongly inner
$H^*$-action. Thus, $A\cong
A^\sigma\otimes\End(H^{\rm op})$ is a matrix algebra, hence the action is inner.
The elements $F_h$ for $h\in H^*$ realizing the action are
subject to the relations: $F_h\cdot F_k=\sum
F_{\h1\k1}\sigma(\h1,\k2)$. One can show that if $f_h$ (respectively, $f'_h$) for $h\in H^*$ realize the action
 on $A^\sigma$ (respectively, $\End(H^*)$), then $F_g=f_{g}\sharp
 f'_{u^{|g|_\sigma}g}$ and $F_v=1\sharp f'_v$. 
The image of
$[A]$ through $\pi\psi$ is, 
as in the previous cases, isomorphic as an algebra to the subalgebra
generated by the $F_h$'s, with opposite product. The
coaction $\rho$ is the restriction of the coaction to this subalgebra; the MUVO action is
given by $a^*\rightharpoonup F_g=\delta_{a, g^{-1}}F_g$ and
$a^*\rightharpoonup F_v=\delta_{a, uv}+\delta_{a,u}F_v$. The $\cal
A$-bicomodule structure behaves on $F_g$, for $g\in G$, as in the case
$V=0$ and $R=R_u$, and it behaves on $F_v$, for $v\in V$, as in the
case $V\neq0$ and $\sigma$ lazy, trivial on $k[G]$. The
assignment $F_h\mapsto F_{S(h)}$ determines an algebra isomorphism 
from $\pi\psi([A])\to k\sharp_{\sigma'}H^*$, where $\sigma'$, the actions and
the coactions are as in the previous cases. In particular, if
$\sigma$ is lazy, i.e., if it is doubly $U$-invariant, then
$|g|_\sigma=0$ for every $g\in G$ and we obtain the same type of formulas that
we obtained for $\sigma$ lazy and trivial on $k[G]$. 
%
%
%
%
%
%

\smallskip

If $G\cong U\times G/U$ then we take the representative
$A=C(1)\sharp\End(H^{\rm op})$ of the class in $BC$, with product
$\bullet$. It corresponds in $BM$, to the representative $A^{\rm op}\cong
(C(1)\sharp\End(H^{\rm op}))^{\rm op}$ which is isomorphic
through $\tau$, as a
module algebra, to the algebra 
$B=\End(H^{\rm op})^{\rm op}\sharp C(1)\cong\End(H^{\rm
op})^{\rm op}\otimes C(1)$. Let us denote by $f_h$, for $h\in H^*$ the
elements realizing the $H^*$-action on $\End(H^{\rm op})^{\rm op}$,
with product $\circ$. 
It is not hard to verify that $A_0= B_0$ is generated by
the elements of the form: $F\sharp 1$ with $F$ commuting with $f_h$ for
every $h\in H^*$, and $F'\sharp x$ with $F'$ commuting with $f_g$ for
$g\in G/U$ and skew-commuting with $f_u$ and $f_v$. Thus,
$\pi\psi([A])$ is generated by the elements
$1\sharp f_{gv_P}$ for $g\in G/U$ and $x\sharp f_u$ in $A$. 
The element $x\sharp f_u$ is central in $A$ because it
corresponds to the element $1\otimes x$ through the above algebra isomorphisms.
By \cite[Lemma 2.3.1 (a)]{CVZ},
$a^*\rightharpoonup(x\sharp f_u)=\delta_{a,1}(x\sharp f_u)$. 
As in the case with $V=0$, for every $g\in G/U$ and every $a\in H^*$
we have
$a^*\rightharpoonup (1\sharp f_g)=\delta_{a,g^{-1}}(1\sharp f_g)$.
The MUVO-action on the elements of the form $1\sharp f_v$ is slightly
more complicated to compute. Let
$\beta(\sum F_i(a)\otimes_{A_0}f_i(a))=1\otimes a^*$, and let
$f_i(a)=1\sharp s_i(a)+x\sharp t_i(a)$. A computation similar to the
previous cases, using that $\rho(x\sharp 1)=x\sharp 1\otimes
\chi$ with $\chi=\sum_{g\in G/U}(g^*-(ug)^*)$, that the $H^*$-action is strongly inner
on the opposite of $\End(H^{\rm op})$, and that if $M$ and $N$ are
$H$-comodules with $\rho_M(m)=\sum \m0\otimes\m1$ and
$\rho_N(n)=\sum\n0\otimes \n1$ then $\rho_{M\sharp N}(m\sharp n)=\sum
\m0\sharp\n0\otimes \n1\m1$ gives 
$$
\begin{array}{rl}
\delta_{a, uv}&=a^*\rightharpoonup(1\sharp
f_v)+\sum F_i(a)\bullet(1\sharp f_u)\bullet(1\sharp
s_i(a))\bullet(1\sharp f_{uv})\\
&-\sum F_i(a)\bullet(1\sharp f_u)\bullet(x\sharp t_i(a))\bullet
(1\sharp f_{uv}).
\end{array}
$$
A computation similar to the previous ones shows that
$$
\begin{array}{rl}
\delta_{a, u}&=\sum F_i(a)\bullet f_{i0}(a)\langle f_{i1}(a),
u\rangle\\
&=\sum F_i(a)\bullet(1\sharp f_u)\bullet(1\sharp
s_i(a))\bullet(1\sharp f_{u})\\
&-\sum F_i(a)\bullet(1\sharp f_u)\bullet(x\sharp
t_i(a))\bullet(1\sharp f_{u})
\end{array}
$$
so, since $(1\sharp f_{uv})=1\sharp(f_u\circ f_v)=-(1\sharp f_u)\bullet
(1\sharp f_v)$, we have
$$a^*\rightharpoonup (1\sharp f_v)=\delta_{a,uv}+\delta_{a,u}(1\sharp f_v).$$
The actions $-\!\triangleright$ and $\triangleleft\!-$ are
as follows:
$$a^*-\!\triangleright\,(x\sharp f_u)=\delta_{u,a}(x\sharp f_u)=(x\sharp f_u)\,\triangleleft\!-\,a^*;$$
$$
\begin{array}{rl}
a^*-\!\triangleright\,(1\sharp f_g)=\delta_{a,g^{-1}}(1\sharp
f_u)=(1\sharp f_u)\,\triangleleft\!-\,a^*;
\end{array}
$$
where we used that $\Delta(a^*)$ is a linear combination of terms of
the form $r^*\otimes s^*$ with $rs=a$; that $S(a^*)$ is a linear
combination of terms of the form $r^*$ with $r\in k[G]\ltimes
\wedge^tV$ if $a\in k[G]\ltimes
\wedge^tV$; that if $a\in k[G]$ the formulas are as in the case
of $V=0$. 
Similarly, for elements in $V$  we have:
$$
\begin{array}{rl}
a^*-\!\triangleright\,(1\sharp f_v)=
\delta_{a,uv}+\delta_{a,1}(1\sharp
f_v)=(1\sharp f_v)\,\triangleleft\!-\,a^*.
\end{array}
$$
The assignment $(x\sharp f_u)\mapsto u$; $(1\sharp f_h)\mapsto
h^{-1}$ and $(1\sharp f_v)\mapsto uv$ determines an algebra
isomorphism $\xi\colon (k[U]\otimes(k[G/U]\ltimes\wedge V))^{\rm
  op}\to k[U]\otimes(k[G/U]\ltimes\wedge V)$. Then
$\pi\psi([C(1)])\cong k[G]\ltimes\wedge (k_\chi\otimes V)$ 
with $k_\chi$ the $1$-dimensional $G$-module corresponding to the
character $\chi$, generalizing \cite[Theorem 5.7, ii]{sequence}.
The Galois object structure is determined by the $H$-coaction: 
$$\rho(u)=u\otimes\chi;\quad\rho(gv_P)=\sum_{a\in
  H^*}\left(\sum\a2(gv_P)S^{-1}(\a1)\right)\otimes a^*;$$
the $H^*$-action:
$$
\begin{array}{l}
a^*\rightharpoonup g=\delta_{a,g}g;\quad a^*\rightharpoonup u=\delta_{a,1}u;\\
a^*\rightharpoonup uv=\delta_{a,uv}+\delta_{a,u}uv;
\end{array}
$$
for every $a^*\in H$, $g\in G/U$, and  $v\in V$; and the $\cal A$-coactions:
$$\rho_1(v)=1\otimes v+v\otimes 1,\quad\rho_2(v)=1\otimes v+v\otimes 1;$$ 
$$\rho_1(g)=g\otimes g,\quad\rho_2(g)=g\otimes g,\quad \forall g\in G.$$

\section{Example: Weyl groups}

In this Section we shall explicitely compute $H^2_L(H)$ and $BM(k, H,
R)$ for a particular family of modified supergroup algebras. We shall
assume that $k={\mathbb C}$. Let $\Phi$ be an
irreducible root system and let $W(\Phi)$ be the corresponding 
Weyl group. Let $V={\mathfrak h}^*$ be its (complexified) natural
representation (or 
its dual), obtained extending the action on $\Phi$ by linearity.  
Let $\{\alpha_1,\ldots\alpha_n\}$ be a fixed fundamental system of simple roots and let $s_1,\ldots,
s_n$ be the corresponding reflections generating $W(\Phi)$ 
as a Coxeter group. For further details the reader is referred to
\cite{bou} and \cite{hu2}.

It is well-known that $V$ is faithful and
irreducible (\cite[Corollary 5.5, Corollary 6.2, Proposition
6.3]{hu2}). By \cite[Proposition 6.3, Lemma and Theorem 6.4]{hu2} the
non-zero,
real, symmetric, $W(\Phi)$-invariant matrices are all equal up to a
scalar. Since the matrices of the representation corresponding to the
basis $\{\alpha_1,\ldots\alpha_n\}$ are real, the same property holds for 
non-zero,
complex, symmetric, $W(\Phi)$-invariant matrices and we have 
$S^2(V^*)^{W(\Phi)}\cong {\mathbb C}$.

Let $w_0$ be the longest element in $W(\Phi)$ (\cite[Section
  1.8]{hu2}).  If $\Phi{\rm {\rm }}$ is of type 
$A_1$, $B_n$, $D_n$ ($n$ even), $E_7, E_8,
F_4, G_2$ (\cite[Corollary 3.19]{hu2}), then $w_0$ acts as $-1$ and
we shall put $G(\Phi)=W(\Phi)$ and $u=w_0$. In this case, 
by \cite[Corollary 3.19]{hu2} for any
ordering $i_1,\ldots, i_n$ of $\{1,\ldots, n\}$ we have
$u=w_0=(s_{i_1}\cdots s_{i_n})^{\frac{\tt
    h}{2}}$ where $\tt h$ is the Coxeter number of $W(\Phi)$ (see \cite[Table on page 80]{hu2}). The
    presentation of $W(\Phi)$ (\cite[Chapitre 4.1.3, Chapitre VI.4.1
    Th\'eor\`eme 1]{bou}) shows that $U=\langle u\rangle$ is a direct summand of
    $G(\Phi)$ if and only if there exists a group morphism $\chi\colon
    G(\Phi)\to k$ with $\chi(u)=-1$. This happens exactly 
when $\Phi$ is of type $A_1$, $B_n$ for $n$ odd, $E_7$ and
    $G_2$.   

If $\Phi$ is of type $A_n$ ($n\ge2$),
$D_n$ ($n$ odd) or $E_6$, then $w_0$ does not act as $-1$ on $V$. 
In this case there exists 
an automorphism $\vartheta$ of order $2$ of the Dynkin diagram
corresponding to $\Phi$ such that, if we extend linearly the action of
$\vartheta$ on $V={\rm span}\{\alpha_1,\ldots\alpha_n\}$, 
the composition of the actions of $\vartheta$
and $w_0$ is $-1$. The subgroup ${\rm
  Aut}(\Phi)$ of ${\rm GL}(V)$ leaving $\Phi$ invariant is the
semi-direct product of its normal subgroup $W(\Phi)$ and the subgroup
of automorphisms of the Dynkin diagram (\cite[\S III.9.2, \S III.12.2]{hu}). So, when $w_0\not=-1$ on $V$
we define 
$G(\Phi)$ as the subgroup of ${\rm Aut}(\Phi)$ 
generated by $\vartheta$
and $W(\Phi)$. Then $G(\Phi)$ is the semi-direct product of $W(\Phi)$ and
$\langle\vartheta\rangle$. 
The representation  $V$, viewed as a $G(\Phi)$-representation 
is still irreducible and faithful. 
For if $\vartheta w=\vartheta w_0w_0 w$ acted trivially on $V$
then $w_0w$ would act as $-1$. Since there is no element in $W(\Phi)$ different
from $w_0$ mapping all positive roots to negative roots, this is
possible only when $w=1$ and $w_0$ acts as $-1$, which is not the case. 
The structure of $G(\Phi)$ is determined by
the action of $\vartheta$ on $W(\Phi)$
that is $\vartheta. w=\vartheta w_0 w_0ww_0w_0\vartheta=w_0ww_0$. In particular, $\vartheta$ and $w_0$
commute and $\vartheta w_0$ is a central involution which we shall
denote by $u$. 
Since $G(\Phi)$ is also generated by the normal subgroup $W(\Phi)$ and
the central subgroup $U=\langle u\rangle$, we have $G(\Phi)\cong
W(\Phi)\times U$. 
%
%
%
%

\smallbreak

Let $H(\Phi):=k[G(\Phi)]\ltimes\wedge V$ with $G(\Phi), V, u $ as
before. We recall that in the $A_1$ case we
obtain Sweedler's Hopf algebra $H_4$. By Theorem \ref{pigra} we have:
%
%
%
$$H^2_L(H(\Phi))\cong H^2(G(\Phi)/U,{\mathbb C}^\cdot)\times {\mathbb C}$$
that is
$$
H^2_L(H(\Phi))\cong
\cases{H^2(W(\Phi)/U,{\mathbb C}^\cdot)\times
{\mathbb C}& for 
  $A_1;B_{n}, n\ge2; D_{2m}, m\ge2;$\cr
& and for $E_7;E_8; F_4; G_2;$\cr
H^2(W(\Phi), {\mathbb C}^\cdot)\times {\mathbb C}
& for $A_{n}, n\ge2;D_{2m+1}, m\ge2;E_6$.\cr}
$$

We shall compute these groups explicitely. The Schur multipliers for
irreducible Coxeter groups have been determined in \cite{IY}. For
the Weyl groups we have:
$$
H^2(W(\Phi),{\mathbb C}^\cdot)\cong\cases{1& for $A_1;A_2;$\cr
{\mathbb Z}_2& for $A_n, n\ge3;B_2;E_6;E_7;E_8;G_2;$\cr
{\mathbb Z}_2\times{\mathbb Z}_2 & for $B_3;D_n, n\ge5; F_4;$\cr
{\mathbb Z}_2\times{\mathbb Z}_2\times{\mathbb Z}_2 & for $B_n, n\ge4;
D_4$
\cr}
$$
so we know $H^2_L(H(\Phi))$ when $w_0\neq-1$ on $V$, i.e., when
$G(\Phi)\cong W(\Phi)\times U$. 

\smallskip

Let $G(\Phi)=W(\Phi)$. If $\Phi=A_1$ then $W(\Phi)=U$ and we have
nothing to prove. If $\Phi$ is of type $B_n$ for odd $n$,
$E_7$ or $G_2$ then  $U$ is a direct summand of $W(\Phi)$ and 
formula (\ref{SY})  
allows us to deduce  $H^2(W(\Phi)/U,{\mathbb C}^\cdot)$  
from the knowledge of $H^2(W(\Phi),{\mathbb C}^\cdot)$ and the
analysis of $\Hom(W(\Phi)/U,{\mathbb Z}_2)$. In these cases we always have
$\Hom(W(\Phi),{\mathbb Z}_2)/\Hom(W(\Phi)/U,{\mathbb Z}_2)
\cong{\mathbb Z}_2$ with
$\Hom(W(\Phi)/U,{\mathbb
Z}_2)\cong{\mathbb Z}_2$ for $\Phi=B_{2m+1}; G_2$ and
$\Hom(W(E_7)/U,{\mathbb
Z}_2)\cong 1$.

\smallskip

If $\Phi$ is of type $B_n$ for $n$ even, $D_n$ for $n$ even, $E_8$ or
$F_4$ then $U$ is not a direct summand of $G(\Phi)=W(\Phi)$ but
sequence 
(\ref{read}) yields:
\begin{equation}\label{quotient}
\CD
1\longrightarrow{\mathbb Z}_2\longrightarrow
 H^2(W(\Phi)/U,{\mathbb
C}^\cdot)\longrightarrow H^2(W(\Phi), {\mathbb C}^\cdot)\longrightarrow
\Hom(W(\Phi),{\mathbb Z}_2).
\endCD
\end{equation} 

Let $\Phi=B_2$. The exact sequence (\ref{quotient}) becomes:
$$
\CD
1\longrightarrow{\mathbb Z}_2@>{\tt T}>> H^2(W(\Phi)/U,{\mathbb
C}^\cdot)@>{\rm Infl}>>{\mathbb Z}_2@>\theta>>{\mathbb
Z}_2\times{\mathbb Z}_2.
\endCD
$$
We would like to describe the image of $\theta$. Let $P$ denote a
lift of a projective representation corresponding to a cocycle $\sigma$ and let
$T_i:=P(s_i)$ for every $i$. Then
$P(w_0)$ is a scalar multiple of $T_1T_2T_1T_2$ and 
centrality of $w_0$ ensures that 
$$P(w_0)T_i=T_iP(w_0)\sigma(s_i,w_0)\sigma^{-1}(w_0,s_i)=
T_iP(w_0)\theta(\sigma)(s_i,u).$$ 
The relations among the $T_i$'s are described in \cite[Table
7.1]{kar}. If $\sigma$ represents the non-trivial class in
$H^2(W(\Phi), {\mathbb C}^\cdot)$ we find $\theta(\sigma)(s_i,u)=-1$ 
for $i=1,2$. It follows that $\Ker(\theta)={\rm
Im}({\rm Infl})=1$ so $H^2(W(\Phi)/U,{\mathbb
C}^\cdot)\cong {\mathbb Z}_2$. 

\smallbreak

Let $\Phi$ be $B_{2m}$ for $m\ge2$. The exact sequence
 (\ref{quotient}) becomes:
$$
\CD
1\to{\mathbb Z}_2@>{\tt T}>> H^2(W(B_{2m})/U,{\mathbb
C}^\cdot)@>{\rm Infl}>>{\mathbb Z}_2\times{\mathbb Z}_2\times{\mathbb
  Z}_2@>\theta>>{\mathbb 
Z}_2\times{\mathbb Z}_2.
\endCD
$$
As before, we would like to understand the map $\theta$. Let $P$
denote the lift of a
projective representation corresponding to a cocycle $\sigma$ and let
$T_i:=P(s_i)$ for every $i$. Then
$P(w_0)$ is a scalar multiple of $(T_1\cdots
T_{2m})^{2m}$ and we have 
$P(w_0)T_i=T_iP(w_0)\theta(\sigma)(s_i,u).$
By the relations in \cite[Table
7.1]{kar} we see that $\theta$ maps the cocycle corresponding to 
$(a,b,c)$ with $a,b,c\in\{\pm1\}$ to the morphism $W(\Phi)\to{\mathbb
  Z}_2$ mapping each generator to $c$. In other words, the only
non-trivial element in ${\rm Im}(\theta)$ is
$\epsilon(w)=(-1)^{\ell(w)}$ where $\ell$ denotes the length function
on a Coxeter group. 
Hence $H^2(W(B_{2m})/U,{\mathbb
C}^\cdot)/{\mathbb Z}_2\cong{\mathbb Z}_2 \times{\mathbb
  Z}_2$. We would like to find out whether  $H^2(W(B_{2m})/U,{\mathbb
C}^\cdot)\cong{\mathbb Z}_2\times{\mathbb Z}_2 \times{\mathbb
  Z}_2$ or $H^2(W(B_{2m})/U,{\mathbb
C}^\cdot)\cong{\mathbb Z}_2\times{\mathbb Z}_4$.  

As in the proof of \cite[Lemma 7.2.8, Lemma 7.2.9]{kar} if
$\sigma$ is a cocycle for $W(B_{2m})/U$ and if $T_i$ denotes the image
of $s_i$ for a projective representation corresponding to $\sigma$,
then the following relations hold, with $n=2m$:
$$
\begin{array}{ll}
T_i^2=1\quad 1\le i\le n,&(T_iT_{i+1})^3=1 \quad 1\le i\le n-2,\\
(T_iT_j)^2=a\quad 1\le i<j+1\le n,&(T_iT_n)^2=b\quad 1\le i\le n-2,\\
(T_{n-1}T_n)^4=c,&(T_1T_3\cdots T_{n-1}T_2\cdots T_n)^n=\lambda\\
\end{array}
$$
where $a,\,b,\,c\in\{\pm1\}$ and $\lambda\in{\mathbb C}^\cdot$. 
We know that $T_i(T_1\cdots T_n)^n=c(T_1\cdots T_n)^nT_i$ by the
previous computation. Since $(T_1\cdots T_n)^n$ and
$(T_1T_3\cdots T_{n-1}T_2\cdots T_n)^n$ are both lifts of $w_0$, they
must differ only by a sign because the relations among the $T_i$'s are the
Coxeter ones up to a sign. Thus, necessarily $c=1$.
Besides, since $u^2=w_0^2=1$ in $W(B_{2m})$, then for its lift there
must hold: 
$(T_1T_3\cdots T_{n-1}T_2\cdots T_n)^{2n}=\lambda^2=\pm1$. If $\psi$ and
$\phi$ are two distinct homomorphisms $W(B_{2m})\to{\mathbb Z}_2=\{0,1\}$ then
$\sigma(wU, w'U)=(-1)^{\phi(w)\psi(w')}$ is a $2$-cocycle on
$W(B_{2m})/U$ corresponding to $(a,b,\lambda)=(1,-1,1)$, as one can
check computing $\sigma(r,s)\sigma^{-1}(s,r)$ as in \cite[Lemma 7.2.8]{kar}. 
The standard representation for $W(B_{2m})$ is a projective
representation for $W(B_{2m})/U$ corresponding to
$(a,b,\lambda)=(1,1,-1)$. If we prove that a triple
$(-1,1,\pm1)$ is also represented, then $H^2(W(B_{2m})/U,{\mathbb
C}^\cdot)\cong{\mathbb Z}_2\times{\mathbb Z}_2 \times{\mathbb
  Z}_2$ because it has four elements of order $2$. By \cite[\S
  7]{stem}, up to a different notation, the representations
corresponding to $(a,b,c)=(-1,1,1)$ are those with $T_n=\pm1$ and
correspond to the projective, nonlinear representations of the
symmetric group $S_{2m}$, i.e., linear representations of its double
cover ${\tilde S}_{2m}$ generated by the $s_i$'s and $-1$. 
This means that the relations become, with $n=2m$:
$$
\begin{array}{l}
T_i^2=1,\quad 1\le i\le n-1\quad(T_iT_{i+1})^3=1 \quad 1\le i\le n-2;\quad
T_n=\pm1\\
(T_iT_j)^2=-1,\quad 1\le i<j+1\le n;\quad(T_1T_3\cdots T_{n-1}T_2\cdots
T_{n-2})^n=\lambda.\\
\end{array}
$$
It is well-known that $w_0'=(s_1s_3\cdots s_{n-1}s_2s_4\cdots
s_{n-2})^m$ is the longest element of $W(A_{n-1})=S_{2m}$. Since $w_0'$ is an involution and since the relations for
the $T_i$'s are the Coxeter ones up to a sign we have
$(T_1T_3\cdots T_{n-1}T_2T_4\cdots T_{n-2})^{n}=\pm1$ 
and the
statement. 
%

\smallbreak

Let $\Phi$ be $D_{4}$. The exact sequence (\ref{quotient}) becomes:
$$
\CD
1\longrightarrow{\mathbb Z}_2@>{\tt T}>> H^2(W(D_4)/U,{\mathbb
C}^\cdot)@>\Infl>>{\mathbb Z}_2\times{\mathbb Z}_2\times{\mathbb Z}_2@>\theta>>{\mathbb
Z}_2.
\endCD
$$
In this case $\theta$ is surjective and it 
maps $(a,b,c)$ with $a,b,c\in\{\pm1\}$ to $abc$. Therefore
$H^2(W(D_4)/U,{\mathbb 
C}^\cdot)/{\mathbb Z}_2\cong{\mathbb Z}_2\times{\mathbb Z}_2$.
We would like to find out whether  $H^2(W(D_{4})/U,{\mathbb
C}^\cdot)\cong{\mathbb Z}_2\times{\mathbb Z}_2 \times{\mathbb
  Z}_2$ or $H^2(W(D_{4})/U,{\mathbb
C}^\cdot)\cong{\mathbb Z}_2\times{\mathbb Z}_4$.  

As in the previous case if
$\sigma$ is a cocycle for $W(D_4)/U$ and if $T_i$ denotes the image
of $s_i$ for the lift of a projective representation corresponding to $\sigma$,
then the following relations hold:
$$
\begin{array}{llll}
T_i^2=1&1\le i\le 4&(T_iT_2)^3=1 &i\neq2\\
(T_1T_3)^2=a& (T_1T_4)=b&(T_3T_4)^2=c,&(T_1T_3T_4T_2)^3=\lambda\\
\end{array}
$$
where $a,\,b,\,c\in\{\pm1\}$ and $\lambda\in{\mathbb C}^\cdot$. 
Since $T_i(T_1T_3T_4T_2)^3=abc(T_1T_3T_4T_2)^3T_i$ we have $c=ab$. 
A direct computation using the relations 
shows that $(T_1T_3T_4T_2)^6=1$, hence $\lambda=\pm1$, and 
$H^2(W(D_{4})/U,{\mathbb
C}^\cdot)\cong{\mathbb Z}_2\times{\mathbb Z}_2 \times{\mathbb
  Z}_2$. Here again, the image of the transgression map is achieved by
the standard representation of $W(D_4)$. 
%

\smallskip 

Let $\Phi$ be $D_{2m}$ for $m\ge3$.  The exact sequence
(\ref{quotient}) becomes: 
$$
\CD
1\longrightarrow{\mathbb Z}_2@>{\tt T}>> H^2(W(D_{2m})/U,{\mathbb
C}^\cdot)@>\Infl>>{\mathbb Z}_2\times{\mathbb Z}_2@>\theta>>{\mathbb
Z}_2.
\endCD
$$
Using the relations in \cite[Table 7.1]{kar} a direct computation
shows that $\theta$ maps $(a,b)$ to $b$ so 
$H^2(W(D_{2m})/U,{\mathbb 
C}^\cdot)/{\mathbb Z}_2\cong{\mathbb Z}_2$. We shall see that, in
analogy with the previous cases $H^2(W(D_{2m})/U,{\mathbb 
C}^\cdot)\cong{\mathbb Z}_2\times{\mathbb Z}_2$. If
$\sigma$ is a cocycle for $W(D_{2m})/U$ and if $t_i$ denotes the image
of $s_i$ for the lift of a projective representation corresponding to $\sigma$,
then the following relations hold, with $n=2m$:
$$
\begin{array}{llll}
t_i^2=1&1\le i\le n;&(t_it_{i+1})^3=1 &1\le i\le n-1;\\
(t_{n-2}t_n)^3=1,&(t_it_j)^2=a&1\le i<j\le n,&i\neq1;\\
(t_{n-1}t_n)^2=b,&&(t_1\cdots t_n)^{n-1}=\lambda;\\
\end{array}
$$
where $a,\,b\in\{\pm1\}$ and $\lambda\in{\mathbb C}^\cdot$. 
Besides, since $t_i(t_1\cdots t_n)^{n-1}=b(t_1\cdots t_n)^{n-1}t_i$ we
necessarily have $b=1$. By \cite[Appendix]{stem} all projective
representations of $W(D_{2m})$ are restrictions of projective
representations of $W(B_{2m})$ with $t_i=T_i$ for $1\le i\le n-1$ and 
$t_n=T_nT_{n-1}T_n$. 
In particular, the representation
corresponding to $(a,b)=(-1,1)$ can be built as the restriction of a 
representation of
$W(B_{2m})$ corresponding to $(a,b,c)=(-1,1,1)$. We have seen
that such representations come from projective representations of
$S_{2m}$ by letting $T_n=\pm 1$, that is, $t_n=t_{n-1}$. We will have our
statement if we show 
that for $(a,b)=(-1,1)$ we can only have $\lambda=\pm1$. We have 
$\lambda=(t_1\cdots t_n)^{n-1}=(t_1\cdots t_{n-2})^{n-1}$. It is not
hard to verify that in $S_{n-1}$ the element $(s_1\cdots
s_{n-2})^{n-1}=1$ so an argument similar to the one used in the
previous cases shows that $\lambda=\pm1$.  

%
%
\smallskip

Let $\Phi$ be $F_4$. The exact sequence (\ref{quotient}) becomes:
$$
\CD
1\longrightarrow{\mathbb Z}_2@>{\tt T}>> H^2(W(F_4)/U,{\mathbb
C}^\cdot)@>\Infl>>{\mathbb Z}_2\times{\mathbb Z}_2@>\theta>>{\mathbb
Z}_2\times{\mathbb Z}_2.
\endCD
$$
The relations in \cite[Table 7.1]{kar} show that $\theta$ maps the
pair $(a,b)$ of 
${\mathbb Z}_2\times{\mathbb Z}_2$ to $(b,b)$, that is, the only
non-trivial element in the image of $\theta$ is, as before, 
$\epsilon(w)=(-1)^{\ell(w)}$. Hence ${\rm
Im}(\Infl)\cong{\mathbb Z}_2$ and we have an exact sequence
$$
\CD
1\longrightarrow{\mathbb Z}_2\longrightarrow H^2(W(F_4)/U,{\mathbb 
C}^\cdot)@>\Infl>>{\mathbb Z}_2\longrightarrow 1.
\endCD
$$ 
Let $\chi\colon
W(F_4)\to{\mathbb Z}_2=\{0,1\}$ be the morphism mapping $s_1$ and
$s_2$ to $1$ and $s_3$ and $s_4$ to $0$.
The map $c(x,y)=(-1)^{\ell(x)\chi(y)}$ is a
  $2$-cocycle on $W(F_4)$. Since $\Hom(W(F_4)/U,{\mathbb
  C}^\cdot)\cong\Hom(W(F_4),{\mathbb C}^\cdot)$, the cocycle $c$ is
also a cocycle for $W(F_4)/U$. Besides, a direct computation using
\cite[Table 7.1]{kar} shows that $\Infl(\bar c)$ is non-trivial because it
corresponds to the pair $(-1,1)$. Since $c^2=1$, the class $\bar c$ determines a section
${\mathbb Z}_2\to H^2(W(F_4)/U,{\mathbb 
C}^\cdot)$ splitting $\Infl$.  

\smallbreak
 
Let $\Phi$ be $E_8$. The exact sequence (\ref{quotient}) becomes:
$$
\CD
1\longrightarrow{\mathbb Z}_2@>{\tt T}>> H^2(W(E_8)/U,{\mathbb
C}^\cdot)@>\Infl>>{\mathbb Z}_2@>\theta>>{\mathbb
Z}_2.
\endCD
$$
Then $\theta(a)=a$,
that is, the non-trivial element in $H^2(W(E_8),{\mathbb C}^\cdot)$
maps to
$\epsilon(w)=(-1)^{\ell(w)}$. Hence $\theta$ is injective, $\Infl$ is
trivial
 and $H^2(W(E_8)/U,{\mathbb 
C}^\cdot)\cong{\mathbb Z}_2$. 

\smallbreak
 
\noindent We have reached the following result:
$$
H^2_L(H(\Phi))\cong\cases{{\mathbb C}& for $A_1; A_2; G_2;$\cr
{\mathbb Z}_2\times {\mathbb C}& for $A_n, n\ge3; B_2; B_3; E_6; E_7; E_8;$\cr
{\mathbb Z}_2\times{\mathbb Z}_2\times{\mathbb C}&  for
$B_{2m+1}, m\ge2; D_n, n\ge5; F_4;$\cr
{\mathbb Z}_2\times{\mathbb Z}_2\times{\mathbb Z}_2\times{\mathbb C}&
for $B_{2m}, m\ge2; D_4.$
\cr}
$$
Next we would like to determine $BM({\mathbb C}, H, R_u)$. By Theorem
\ref{tutte} and Proposition \ref{coincide} we have:
$$BM({\mathbb C}, H(\Phi), R_u)\cong
\cases{{\mathbb Z}_2\times H^2(W(\Phi), 
    {\mathbb C}^\cdot)\times {\mathbb C}& for $A_1;B_{2m+1};E_7;
    G_2$\cr
{\mathbb Z}_2\times H^2(G(\Phi), 
    {\mathbb C}^\cdot)\times {\mathbb C}& for $A_{n}, n\ge2;D_{2m+1};E_6$\cr
H^2(W(\Phi), {\mathbb C}^\cdot)\times{\mathbb C} & for
$B_{2m};D_{2m};E_8; F_4$.\cr}
$$
%
%

Let $\Phi=A_n, n\ge2;D_n, n\ {\rm odd}, E_6$. Then  $G(\Phi)\cong
W(\Phi)\times U$ and (\ref{SY}) yields
$$H^2(G(\Phi), k^\cdot)\cong H^2(W(\Phi), k^\cdot)\times {\mathbb Z}_2.$$

We can conclude with:
$$
BM({\mathbb C}, H, R_u)\cong
\cases{{\mathbb Z}_2\times {\mathbb C}& for $A_1; A_2; B_{2};E_8;$\cr
{\mathbb Z}_2\times{\mathbb Z}_2\times {\mathbb C}& for $A_n, n\ge3;
D_{2m}, m\ge 3;$ \cr
& and $E_6; E_7; F_4;G_2;$\cr  
{\mathbb Z}_2\times{\mathbb Z}_2\times{\mathbb Z}_2 \times
{\mathbb C}& for $B_3; B_{2m}, m\ge2;$\cr
& and for $D_4;D_{2m+1}, m\ge2;$\cr
{\mathbb Z}_2\times {\mathbb Z}_2\times{\mathbb Z}_2\times{\mathbb
Z}_2\times {\mathbb C}& for $B_{2m+1},m\ge2.$\cr }
$$

\section*{Acknowledgements}

The author wishes to thank J. Bichon and J. Cuadra for useful discussions.
The research leading to this manuscript has been partially supported by Progetto
Giovani Ricercatori CPDG031245 of the University of Padua.

\end{document}